\documentclass{siamart190516}
\usepackage{braket,amsfonts}
\usepackage{amsmath,amssymb,mathtools,latexsym}
\usepackage{graphicx,epstopdf}
\usepackage{amsopn}
\usepackage{mathrsfs}

\Crefname{ALC@unique}{Line}{Lines} 

\usepackage{algorithm}
\usepackage[noend]{algpseudocode}

\usepackage{array}     
\usepackage{url,tikz,pgfplots,xcolor}
\usepackage[export]{adjustbox}

\DeclarePairedDelimiter\norm{\lVert}{\rVert}%

\usepackage[font={small,singlespacing},skip=1pt]{subcaption}
\usepackage[font={small,singlespacing},skip=1pt]{caption}
\usepackage{afterpage}

\makeatletter
\providecommand{\l@section}{}
\providecommand{\l@subsection}{}
\providecommand{\l@subsubsection}{}
\providecommand{\l@paragraph}{}
\providecommand{\l@subparagraph}{}
\providecommand{\l@figure}{}
\providecommand{\l@table}{}
\usepackage[subfigure]{tocloft}
\makeatother

\usepackage{physics}

\title{Hierarchical Matrix Approximations of Hessians \\ Arising in Inverse Problems Governed by PDEs\thanks{\funding{This work was supported by the King Abdullah University of Science and Technology (KAUST) Office of Sponsored Research (OSR) under Award No: OSR-2018-CARF-3666.}}}

\author{
Ilona \negthinspace Ambartsumyan\thanks{Oden Institute for
  Computational Engineering and Sciences, The University of Texas at Austin. (\email{ailona@austin.utexas.edu}, \email{tanbui@ices.utexas.edu}, \email{omar@ices.utexas.edu}).}
\and \negthinspace Wajih \negthinspace Boukaram\thanks{Extreme Computing Research Center, King Abdullah University of Science and Technology. 
   (\email{wajihhalim.boukaram@kaust.edu.sa}, \email{stefano.zampini@kaust.edu.sa}, \email{david.keyes@kaust.edu.sa}).}
\and \negthinspace Tan \negthinspace Bui-Thanh\footnotemark[2]
\and \negthinspace Omar \negthinspace Ghattas\footnotemark[2]
\and David Keyes\footnotemark[3]
\and Georg Stadler\thanks{Courant Institute of Mathematical Sciences, New York University. (\email{stadler@cims.nyu.edu}).}
\and George Turkiyyah\thanks{Department of Computer Science, American University of Beirut, Lebanon.
  (\email{gt02@aub.edu.lb}).}
\and Stefano Zampini\footnotemark[3]
}

\headers{Hierarchical Hessians}{Hierarchical Hessians}
\hypersetup{ pdftitle={H Hessians} }

\begin{document}
\maketitle

\begin{abstract}
Hessian operators arising in inverse problems governed by partial
differential equations (PDEs) play a critical role in delivering
efficient, dimension-independent convergence for both Newton solution
of deterministic inverse problems, as well as Markov chain Monte Carlo
sampling of posteriors in the Bayesian setting. These methods require
the ability to repeatedly perform such operations on the Hessian
as multiplication with arbitrary vectors, solving linear systems,
inversion, and (inverse) square root. Unfortunately, the
Hessian is a (formally) dense, implicitly-defined operator that is
intractable to form explicitly for practical inverse problems, 
requiring as many PDE solves as inversion parameters. Low rank
approximations are effective when the data contain limited information
about the parameters, but become prohibitive as the data become more
informative. 
However, the Hessians for many inverse problems arising in practical
applications can be well approximated by matrices that have
hierarchically low rank structure. Hierarchical matrix representations
promise to overcome the high complexity of dense representations and
provide effective data structures and matrix operations that have only
log-linear complexity. In this work, we describe algorithms for
constructing and updating hierarchical matrix approximations of
Hessians, 
and illustrate them on a number of representative inverse problems involving
time-dependent diffusion, advection-dominated transport, frequency
domain acoustic wave propagation, and low frequency Maxwell
equations, demonstrating up to an order of magnitude speedup 
compared to globally low rank approximations.
\end{abstract}


%

\section{Introduction}
\label{sec:intro}


The Hessian operator plays a central role in optimization of systems
governed by partial differential equations (PDEs), also known as {\em
  PDE-constrained optimization}. While the approach proposed here
applies more broadly to other PDE-constrained optimization problems
including optimal control and optimal design, we will focus on an
important class: inverse problems.  The goal of an inverse problem is
to infer model parameters, given observational data, a forward model
or state equation (here in the form of PDEs) mapping parameters to
observables, and any prior information on the parameters. Often the
parameters represent infinite-dimensional fields, such as
heterogeneous coefficients (including material properties),
distributed sources, initial or boundary conditions, or geometry. We
focus here on this infinite dimensional setting, leading to large
scale inverse problems after discretization.

The Hessian operator plays a critical
role in inverse problems. For deterministic
inverse problems, finding the parameters that best fit the data is
typically formulated as a regularized nonlinear least squares
optimization problem. Its objective function(al) consists of two
terms: the data misfit, that is, the squared $\ell^2$-norm of the
difference between the output observables predicted by the PDE for
given model parameters and the observed data; and a regularization
term that ensures stability by penalizing unobserved features of the
parameters, such as rough components. The Hessian is given by the
second variation, with respect to model parameters, of this
regularized data misfit. Minimizing the objective by Newton's method
requires solution of a sequence of linear systems with the Hessian as
its system matrix. Due to its affine invariance and resulting
mesh-independent behavior (e.g., \cite{Heinkenschloss93}), Newton's
method is the gold standard for inverse problem solution, with
demonstrated convergence independent of parameter dimension for
large-scale, complex inverse problems (e.g.,
\cite{EpanomeritakisAkccelikGhattasEtAl08, IsaacPetraStadlerEtAl15}).
Efficient solution of the linear systems arising at each Newton
iteration requires an effective preconditioner
that is inexpensive to form and ``invert.'' The challenge is
that the discretized Hessian of the data misfit functional is formally
a dense operator that cannot be formed explicitly, since each column
requires the solution of a linearized state PDE. Instead, we can
exploit the fact that for ill-posed inverse problems, the data misfit
Hessian operator
is compact (its eigenvalues accumulate at zero), and when the
regularization operator takes the typical form of an elliptic
differential operator, we can precondition by the inverse of the
regularization to yield a preconditioned Hessian in the form of a
compact perturbation of the identity. A conjugate gradient linear
solver will converge in a number of iterations that depends on the
number of distinct eigenvalues, which for regularization
preconditioning is small when the compact part has eigenvalues that
decay rapidly, so that it has small effective rank $r$. Since at each
iteration, applying the Hessian to a vector requires the solution of a
pair of linearized state/adjoint PDEs, $O(r)$ such pairs must be
solved at each Newton iteration. As demonstrated in a number of
geophysical inversion problems involving global seismology
\cite{Bui-ThanhGhattasMartinEtAl13}, ice sheet flow
\cite{IsaacPetraStadlerEtAl15}, atmospheric transport
\cite{FlathWilcoxAkcelikEtAl11}, poroelastic subsurface flow
\cite{HesseStadler14},
and joint inversion \cite{CrestelStadlerGhattas18}, $r$ is small and
dimension-independent, and such a
preconditioning strategy is effective. However, as we shall see below,
the central challenge is that $r$ grows as the data become more
informative (a desirable situation for inverse problems), and
regularization preconditioning becomes prohibitive.\footnote{The
  augmented Lagrangian KKT preconditioner
  \cite{AlgerVillaBui-ThanhEtAl17} overcomes this difficulty and is
  insensitive to $r$, but since it involves the optimality system of the
  PDE-constrained optimization problem, is storage-prohibitive for
  large-scale time-dependent inverse problems.}

Another critical role of the Hessian is in Bayesian inversion, which
provides a powerful and systematic framework for quantifying
uncertainties in the solution of inverse problems. Given probability
distributions that represent observational data and their associated
uncertainties, the state PDE model and its associated uncertainty,
and any prior knowledge on the model parameters, Bayesian solution of
the inverse problem yields the posterior distribution, which
characterizes the probability that any particular parameter field gave
rise to the observed data. The Hessian in this setting is given by the
second variation (with respect to model parameters) of the negative
log of the posterior distribution. This is equivalent to the
deterministic Hessian in the common setting of Gaussian additive noise
and prior. For linear inverse problems, the posterior covariance is
equal to the inverse of the Hessian. For moderately nonlinear inverse
problems, the posterior can be approximated by the Laplace
approximation, in which the posterior covariance is given by the
inverse Hessian evaluated at the point that maximizes the posterior
covariance (or minimizes the negative log posterior). Finally, for
highly nonlinear inverse problems, where the Laplace approximation
does not adequately capture the uncertainty in the inverse solution,
Markov chain Monte Carlo (MCMC) methods are used to sample the
posterior, from which a sample posterior covariance can be
constructed. In this case the Hessian plays an instrumental role in
exploiting the geometry of the posterior to accelerate MCMC
convergence \cite{MartinWilcoxBursteddeEtAl12, CuiLawMarzouk16,
  BeskosGirolamiLanEtAl17, Bui-ThanhGirolami14}. In all three cases,
what is needed is an accurate and efficient method to approximate the
Hessian, its inverse, and their square roots.

The contemporary approach to approximation of Hessians of ill-posed
inverse problems is to again exploit the compactness of the Hessian of
the data misfit functional. Compactness can be theoretically proven in
specific settings (e.g., for acoustic and electromagnetic inverse
shape and medium scattering \cite{Bui-ThanhGhattas13a,
  Bui-ThanhGhattas12, Bui-ThanhGhattas12a}), or else demonstrated
numerically for an  array of inverse problems
\cite{BashirWillcoxGhattasEtAl08, MartinWilcoxBursteddeEtAl12,
  Bui-ThanhBursteddeGhattasEtAl12, Bui-ThanhGhattasMartinEtAl13,
  ZhuLiFomelEtAl16, HesseStadler14, PetraMartinStadlerEtAl14,
  IsaacPetraStadlerEtAl15, FlathWilcoxAkcelikEtAl11}.
This often allows one to make a low rank approximation of the data
misfit component of the Hessian (preconditioned by the prior
covariance, which plays the role of the inverse regularization
operator), computed efficiently via matrix-free randomized SVD
\cite{halko11}. This entails $O(r)$ products of the
prior-preconditioned data misfit Hessian matrix with random vectors,
each of which as above requires a pair of linearized state/adjoint
PDE solutions. The inverse Hessian and its square root are then
computed at negligible additional cost using the resulting spectral
decomposition along with the Sherman-Morrison-Woodbury formula
\cite{Bui-ThanhGhattasMartinEtAl13,  PetraMartinStadlerEtAl14}.
When the inverse problem is highly ill-posed, i.e., when the data
inform few components of the parameter field, $r$ is small. Moreover,
the eigenfunctions corresponding to the dominant eigenvalues of the
Hessian are often smooth (as a consequence of the data being unable to
inform rough components of the parameter field), so that $r$ is
independent of the mesh size and hence parameter dimension. In such
cases, the low rank approximation is both efficient and
scalable (with respect to parameter dimension).

In summary, to ensure efficiency of CG solution with regularization
preconditioning as well as low rank approximation of the
prior-preconditioned data misfit Hessian, it is critical that the
effective rank $r$ of the Hessian remain small relative to the
parameter dimension.
However, while $r$ generally does not
grow with increasing parameter dimension, it does grow as the data
become  more informative about the parameters. In fact,
for a linear Bayesian inverse problem with additive Gaussian noise,
the Kullback-Leibler divergence (or relative entropy) from posterior
to prior measure---which quantifies the information gained from the
data---can be shown to be equal to
$ \Sigma_i \log (\lambda_i + 1) $, where $\lambda_i$ are the
eigenvalues of the prior-preconditioned data misfit Hessian \cite{AlexanderianGloorGhattas16}. This sum
can be truncated when $\lambda_i \ll 1$, so the dominant
eigenvalues---and hence effective rank of this operator---are directly
related to the information content of the data.  The information
gained from the data (and hence $r$) generally increases as the
number of sources (or experiments) increases and the number of
observations (or sensors or receivers) increases. The data also become
more informative as the state PDEs become less dissipative (for
example as flow or transport equations become more
advection-dominated), or as they resolve finer scales
(such as in wave propagation with increasing frequency). Finally, as the
noise in the data decreases, the strength of the data misfit Hessian
increases relative to the prior/regularization, again increasing the
effective rank and hence informativeness of the data. In all such cases,
both regularization preconditioning as well as low rank-based data
misfit Hessian approximation become prohibitive, and there
remains a critical need for more efficient Hessian approximation.

In this paper, we show for the first time the efficiency
of hierarchical matrix representations \cite{hackbusch15}
in representing the Hessian operator of various PDE-constrained problems, and
consider these tunable-accuracy approximations to
address the shortcomings of the globally low rank representations
in the increasingly data-informed regime, as well as to provide
log-linear time complexity algorithms for constructing robust approximations
of the inverse of these operators.
Hierarchical matrices exploit the fact that many of their off-diagonal
blocks can be approximated to high accuracy with blocks of low rank,
allowing substantial compression in the memory footprint relative to
the dense representation.  The blocks that may be so represented can
be of different sizes offering a natural tree hierarchy for managing
and operating on them, and resulting in matrices that can be stored
and operated on in linear or log-linear space and time complexity as
opposed to polynomial complexity of dense
representations. Such matrix representations provide an
algebra in which memory versus accuracy tradeoffs may be explicitly made
to suit the needs of specific application contexts, a particularly
useful feature for algorithms that can operate with coarser
approximations to benefit from substantial speedup.  Hierarchical
matrices may be viewed as algebraic generalizations of fast multipole
methods \cite{yokota14}, allowing  the fast performance of not only
matrix-vector products,  but also a full stack of
linear algebra operations including generalizations of randomized
algorithms for constructing matrix decompositions
\cite{halko11,boukaram19}.

With the emergence of GPUs and manycore architectures as key platforms
for high performance scientific computing, efficient execution on
these architectures has become a critical feature for algorithms that
aim to be deployable in practice. With this in mind, the algorithms we
present for operating on hierarchical representations of Hessians are
developed to exploit the high throughput of these modern architectures
through data parallel operations and careful orchestration of data
movement to minimize latencies. Linear algebra operations such
as matrix-vector products, low rank updates, and matrix inversion can be
performed directly on the GPU resident compressed representations,
allowing algorithms to benefit both from the reduced memory footprint
with its resulting log-linear operation count as well as from the high
flop rate of modern hardware architectures.

Besides low rank approximation, several methods for compact
approximation of (data misfit) Hessians stemming from inverse problems
have been developed recently. The pseudodifferential scaling method of
\cite{Nammour11, NammourSymes09} exploits the pseudodifferential
nature of the Hessian in seismic inverse problems, in particular that
it is approximately diagonal in phase space and can thus be estimated
by a single application to a vector. Also exploiting the
pseudodifferential structure of seismic inversion Hessians is the
matrix probing method of \cite{DemanetLetourneauBoumalEtAl12}, which
approximates the Hessian (and its inverse) with basis matrices
stemming from the Hessian's symbol and find their coefficients by
probing the Hessian in random directions. Recently, methods that
exploit the local translation invariance of Hessians have been
introduced \cite{ZhuLiFomelEtAl16, AlgerRaoMeyersEtAl19}. The adaptive
product-convolution approximation in particular is demonstrated to be
robust to the Peclet number for advection-dominated transport and the
frequency for an auxiliary operator that arises in connection with KKT
preconditioning \cite{AlgerRaoMeyersEtAl19} for a wave inverse problem
\cite{Alger19}. Here we focus on comparisons with
low rank approximation, and defer comparison to these other methods
in appropriate contexts for future work.


The rest of this paper is organized as follows. In \cref{sec:why}, we
argue why  Hessians of PDE-governed inverse problems are
expected to be well approximated by hierarchical
matrices. \Cref{sec:Hrep} introduces the $\mathcal{H}^2$ data
structures used to store the hierarchical matrix as well as key
operations on the representation, including matrix-vector
multiplication, compression, low rank updates, and construction from
randomized sampling. 
In \cref{sec:applications} we present numerical
experiments to assess the effectiveness of the hierarchical matrix
Hessian representations on inverse problems governed by diffusion,
transport, and acoustic and electromagnetic wave propagation. The
experiments demonstrate the superiority of hierarchical matrix
approximations over a low rank approximation as the data become more
informative, with speedups ranging from two to over an order of
magnitude. \Cref{sec:conclusions} concludes the paper.


\section{Why do Hessians Admit Hierarchical Low Rank Representations? \nopunct}
\label{sec:why}

The Hessians that arise in PDE-governed inverse problems (whether
deterministic or Bayesian) may be approximated to specified accuracy by
highly compressible matrices, a property they inherit from the
underlying PDE operators \cite{bebendorf03}, often reinforced by
sparse observation operators. As argued above, the global rank of
these Hessians increases as the data become more informative, such
that a low rank approximation may no longer be tenable.
However, as we argue below, these Hessians typically have blocks that
can be well 
approximated by low rank representations with bounded rank. These
blocks are at different levels of granularity because of their
different sizes, making hierarchical matrix representations an
ideal vehicle for storing and computing with Hessians.

In this paper, we consider inverse problems governed by PDEs with
distributed parameter fields to be inferred from data. After
discretization over a $d$-dimensional domain $\Omega$ (which without
loss of generality we assume is conducted by finite elements), we obtain
an optimization problem of the 
form
\begin{equation}
\label{eq:optproblem}          
  \underset{m}{\mbox{minimize }}  {\displaystyle J(m) := 
                        F(u(m)) + \alpha R(m)} 
\end{equation}
where the state variables $u(m) \in \mathbb R^N$ depend on the model
parameters $m \in \mathbb R^n$ via solution of the discretized PDEs
\begin{equation}
                   g(m, u) := K(m) u - f = 0,
\label{eq:pdeconstraint}            
\end{equation}
with coefficient matrix $K \in \mathbb{R}^{N\times N}$ and source $f
\in \mathbb{R}^N$. For simplicity, we confine our discussion to the
case where the forward problem $g(m,u)=0$ represents elliptic PDEs
that have arbitrary dependence on the parameter but are linear in the
state. However, our methodology is more broadly applicable to
parabolic and hyperbolic forward problems, as demonstrated in the
applications in \cref{sec:applications}, and to nonlinear
forward problems. Multiple sources (or multiple experiments) may be
used to generate data, in which case \eqref{eq:pdeconstraint}
is indexed by source $s$, generating the state $u^s$.


The data misfit term $F(u(m))$ measures the difference between recorded
observations and the response of the model at receivers
placed on the boundary or the interior of $\Omega$.
It is often taken to be the sum of squares, over all
sources and receivers, of $d_r^s$, the difference between the observed state
at receiver $r$ due to source $s$, and the
corresponding model response $u^s_r$,
\[ F :=\tfrac{1}{2}\sum_s \sum_r \| u^s_r - d_r^s \|^2 \]
although other measures of misfit may be sometimes preferable
\cite{engquist18}. $R(m)$ is a regularization term expressing prior
information on the parameter field (such as piecewise smoothness) and
$\alpha$ is a scalar regularization parameter that attempts to
annihilate components of $m$ that are uninferable from the data while
preserving those that can be inferred. The regularization term is generally a local
operator whose discretization is sparse and does not pose
computational difficulties even at large scale. Thus in this section
we focus on the treatment of the problematic $F$ term, whose
evaluation requires solution of the forward PDEs for a given $m$ for
each source $s$. 


We define $G \in \mathbb{R}^{N\times n}$ to be the
partial derivative of the forward PDE residual with respect to the
parameter field, i.e.,
\[
G := \partial_m g = \partial_m K \times_2 u,
\]
where $\partial_m K$ is a third order tensor of size $N \times N
\times n$ and $\times_2$ is the 2-mode tensor vector product operation
in the notation of \cite{kolda09}. Since $K$ is a discretization of a
local PDE operator, it is representable as a sparse matrix. The
parameter field is typically discretized on the same mesh as the state
(perhaps at lower order), 
and thus the local operator $G$ is also sparse.

Unlike the local sparse matrices $K$ and $G$, the inverse operator
$K^{-1}$ is a discretized, non-local solution operator, and is in
general formally dense.  It is however data-sparse, since many of its
off-diagonal blocks can be represented to high accuracy by low rank
approximations with bounded rank $k$, and these blocks occur at
different levels of granularity \cite{bebendorf03,borm10}. For a
block $ts$ of $K^{-1}$, with rows and columns corresponding to the 
index set $t$ and $s$ of clusters of nodes of a finite element mesh, the
compressibility condition states that if the bounding boxes of these
clusters are sufficiently far away from each other, then the
local rank of the $ts$ block is essentially independent of the block
size for Poisson-like operators. Such a condition is known as an
admissibility condition, and it depends on the spatial distance
between the clusters of $t$ and $s$, relative to their size.

Using the nested basis $\mathcal{H}^2$ representation \cite{borm07},
the total amount of storage needed to represent $K^{-1}$ to a given
accuracy $\epsilon$ is only $O(Nk)$, where the local block rank $k$
grows as $|\log \epsilon|^{d+1}$ for PDEs in $d$ spatial dimensions.
Note that the local ranks of the blocks are independent of the global
rank.  Even for matrices that are of full rank $N$, the local ranks
grow slowly, with only a logarithmic dependence on target
approximation accuracy, as well as with a logarithmic dependence on
high contrast coefficients \cite{bebendorf16}. High-frequency problems
can be tackled effectively by using directional compression
techniques, with asymptotic storage $O(N k^2 \log N)$ \cite{borm17}.


The gradient of the data misfit term of the objective function may be
written in terms of an adjoint variable $p \in \mathbb{R}^N$, defined
via the adjoint equations \cite{HinzePinnauUlbrichEtAl09}
\begin{equation}
	l(m, p) := K^T\!(m) \, p + \partial_u F = 0.
	\label{eq:adjoint}
\end{equation}
The gradient of $F$ is then given by
\begin{equation}
	\nabla_m F := (\partial_m K \times_2 u)^T p = - G^T K^{-T} \partial_u F.
	\label{eq:gradient}
\end{equation}
We denote by $L \in \mathbb{R}^{N\times n}$ the partial derivative of
the residual of the adjoint equations \eqref{eq:adjoint} with respect
to $m$,
\[
L := \partial_m l = \partial_m K^T \times_2 p,
\]
which is a local operator and has a sparsity pattern similar to that
of $G$. Taking the derivative of \eqref{eq:gradient} with respect to
$m$, we obtain an expression for the Hessian matrix in terms of the state and
adjoint variables $u$ and $p$, and the local operators $G$ and $L$,
\begin{equation}
\nabla^2_{mm} F = G^T K^{-T} \partial_{uu}^2 F K^{-1} G - G^T K^{-T} L
- L^T K^{-1} G + (\partial_{mm}^2 K \times_2 u) \times_1 p.
\label{eq:hessian}
\end{equation}
The first term of the latter expression represents the positive
semi-definite Gauss-Newton part of the Hessian, involving the
application of the local operator $G$ and its adjoint from the left
and right to the triple product $K^{-T} (\partial_{uu}^2 F) K^{-1}$.
The next two terms involve the application of the sparse
Jacobians of the state and adjoint equations, from the left and from
the right, to $K^{-1}$ and its adjoint $K^{-T}$.  The last term
involves a sparse, fourth-order tensor multiplied in 1-mode and 2-mode
by the adjoint and state variables.

Equation \eqref{eq:hessian} shows why it is prohibitive to form
the Hessian operator explicitly: forming $K^{-1} G$ (and similar products) requires a
number of forward PDE solves equal to the number of parameters
multiplied by the number of sources. On the other hand,
\eqref{eq:hessian} also shows that products of the Hessian with
arbitrary vectors can be carried out at the cost of one forward solve
(with $K$) and one adjoint solve (with $K^T$). Note also that (unless
the forward problem is linear in $m$), in
general $u$, $p$, $K$, $G$, and $L$ all depend on $m$, so the Hessian
varies over parameter space. However the argument below is not
predicated on where in parameter space the Hessian is evaluated. 
 
%
%

In the form given in \eqref{eq:hessian}, we can see why the
Hessian of \eqref{eq:optproblem} exhibits a hierarchical low rank
structure.  We first consider the Gauss-Newton term.  The observation
operator $\partial_{uu}^2 F$ is positive semi-definite, and in
particular is block diagonal for the standard weighted least squares
objective (where the blocks correspond to the support of the receiver
operators). When the number of observations in $F$ is small, it could have
relatively small global rank as well.  A general result concerning the
structure of the product of hierarchical matrices \cite[Theorem 2.24]{grasedyck03}
shows that the triple product $K^{-T} \partial_{uu}^2 F K^{-1}$, even
when $\partial_{uu}^2 F$ is full rank, may be represented as a
hierarchical matrix, although the local ranks may grow as a result.
By the same argument, the right and left multiplications of the latter triple product by the
sparse matrix $G$ and its transpose also produce a hierarchical matrix with 
bounded local ranks. For favorable sparsity structures
in $G$, the resulting local ranks may even decrease. 
This is due to the fact that a block in the product
contains contributions from a bounded number of block product pairs, where
most of the operands involved are zero and contribute no data to the
target block.  The local ranks of the Gauss-Newton Hessian may
therefore end up being smaller than the local ranks of the solution
operator.  The second and third terms of the full Hessian expression
given in \eqref{eq:hessian} also have the form of a
hierarchical matrix multiplied on both sides by sparse matrices, and
therefore result in hierarchical matrices with small local
ranks. Finally, the last term of the Hessian involves the sparse
fourth order tensor. The sum of all four terms may of course increase
the local ranks slightly, but the resulting Hessian will still retain
a hierarchical structure.


\section{Hierarchical Matrix Construction of Dense Hessians}
\label{sec:Hrep}


In this section we describe the hierarchical Hessian representation that results in linear space complexity, and show how to perform linear algebra operations on the Hessian directly in its compressed representation, including matrix-vector products and norm computations, in near linear time complexity. We also describe a randomized procedure for constructing hierarchical Hessian approximations to a desired accuracy $\epsilon$.

\subsection{Structure and representation} Various representations for hierarchical matrices have been proposed in the past and may be roughly classified along two characteristics as shown in \cref{tbl:H}. One relates to the structure of the block partitioning of the matrix, and the other relates to the format of the representation and how the low rank data is stored. In the simplest representations, such as $\mathcal{H}_p$ \cite{hackbusch15} or HODLR \cite{darve13} (see also \cite{yu17, yu18}), a weak-admissibility condition is used for partitioning, so that all the off-diagonal blocks touch the main diagonal, and each block at level $l$ of the hierarchy defined by the row and column index sets $t$ and $s$ respectively has its own low rank representation $A_{ts}^l = U_{ts}^l {V_{ts}^l}^T$. The HSS representation \cite{xia10} improves the asymptotic complexity by using nested bases while still keeping the simple partitioning. In nested basis representation, the column and row bases are no longer specific to a given block but shared with all blocks in the $t$ block row and $s$ block column, i.e. $A_{ts}^l = U_{t}^l S_{ts} {V_{s}^l}^T$.  The bases $U_t$ and $V_s$ are not stored explicitly, but are computed recursively on demand from their ``children,'' i.e., block rows and columns that are subsets of $t$ and $s$. Bases are only stored explicitly at the lowest level, and only small transfer matrices are needed to compute coarser levels bases, resulting in the asymptotically optimal $O(kn)$ storage, where $k$ is the maximum rank among $A^l_{ts}$ blocks.

\begin{table}[!ht]
\[
\begin{array}{|l||c|c|}
	\hline
	     & \text{Flat bases} & \text{Nested bases} \\ \hline \hline
	\text{Weak admissibility partitioning} & \text{HODLR}\cite{darve13}, \mathcal{H}_p \cite{hackbusch15} & \text{HSS} \cite{xia10} \\ \hline
	\text{Strong admissibility partitioning} & \mathcal{H} \cite{hackbusch15} & \mathcal{H}^2 \cite{borm07},  \mathcal{DH}^2 \cite{borm17} \\ \hline
\end{array}
\]
\caption{Classification of hierarchical matrix representations.}
\label{tbl:H}
\end{table}

\begin{figure}[!ht]
\begin{center}
	\begin{subfigure}{.66\textwidth}
		\centering
		\includegraphics[width=0.7\linewidth,valign=m]{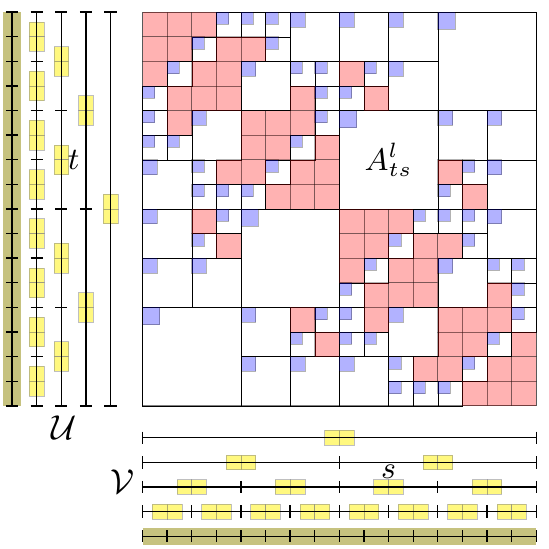}
		\label{fig:h2b}
	\end{subfigure}
	\hspace*{0.05in}
	\begin{subfigure}{.28\textwidth}
		\centering
		\includegraphics[width=\textwidth,valign=m]{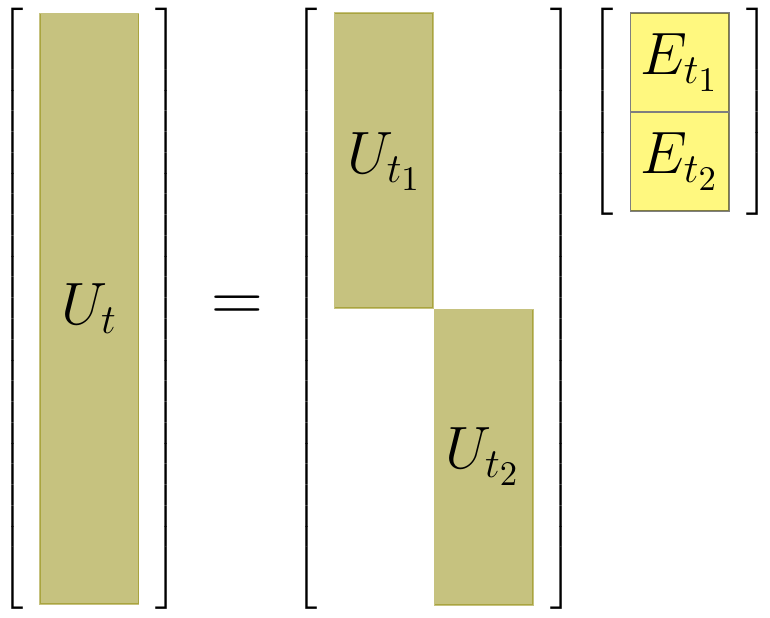}
		\label{fig:h2a}
	\end{subfigure}
\end{center}
\caption{$O(n)$ hierarchical representation of $A$: matrix partitioning and nested bases. Low rank representation of the blocks is of the form $A_{ts}^l = U^l_t S_{ts}^l {V^l_s}^T$, with $U_t^l$ computed from the bases of its children $U^{l-1}_{t1}$ and $U^{l-1}_{t2}$ using small transfer matrices $E^{l-1}$ stored in a column basis tree $\mathcal{U}$. Per-block low rank data $S_{ts}$ are stored as small $k\times k$ blocks in a matrix tree $\mathcal{S}$.}
\label{fig:h2}
\end{figure}

The primary drawback of the weak admissibility
partitioning is that the ranks of the off-diagonal blocks grow too
large even when only a moderate accuracy is requested for many
problems of interest, and in particular for three-dimensional
problems. This undesirable growth in rank demands a more refined
matrix blocking, tuned to the actual geometric discretization of the
problem, to be used for the partitioning of the hierarchical
matrix. In this blocking, identified as strong admissibility partitioning in the bottom row of \cref{tbl:H}, only those matrix blocks $ts$ whose $t$ and $s$ index sets define clusters that are sufficiently far away from each other will admit a low rank representation. This allows more refinements to take place on large blocks of the weak partitioning block structure. In this work, we use the $\mathcal{H}^2$ representation \cite{borm07}; for a schematic representation of the storage format, see \cref{fig:h2}. A general partitioning of the matrix is allowed, as prescribed by an application-dependent admissibility condition, quantifying the relative distance between clusters. Every admissible low rank block is represented by a triple product $U_tS_{ts}V_s^T$: all the needed information for generating the $U_t$ columns bases is stored in the leaves of the $\mathcal{U}$ tree. 
For non-symmetric matrices, a separate row basis tree $\mathcal{V}$ is also needed. Blocks that do not admit a low rank representation are stored as dense blocks, and are limited to the leaf level of the matrix blocks. An informal notation for a hierarchical matrix $A$ is then
\begin{equation}
	A_{\mathcal{H}^2} = D + \mathcal{U} \cdot \mathcal{S} \cdot \mathcal{V}^T
\end{equation}
where $D$ is a block sparse matrix, and the tree triplet
$\{\mathcal{U}, \mathcal{S}, \mathcal{V}\}$ provides the data for the
low rank blocks at multiple levels of granularity. Hierarchical
matrices may hence be viewed as generalizations of matrix
decompositions of the form ``diagonal plus low rank'' or ``sparse plus
low rank'' that have been used with success in PDE-governed inverse
problems \cite{Bui-ThanhGhattasMartinEtAl13}.

\subsection{Fast operations on hierarchical matrices}

The power of hierarchical matrices is not only a result of their optimal storage but of the fact that we can also perform general linear algebra operations directly in the compressed representation in linear or log-linear time complexity.

Matrix-vector multiplication, for example, can be done in four stages:
\begin{equation}
	A_{\mathcal{H}^2} x  = D x + \mathcal{U} \cdot \mathcal{S} \cdot \mathcal{V}^T \cdot x = D x + \mathcal{U} \cdot ( \mathcal{S}  \cdot (\mathcal{V}^T \cdot x)))
\end{equation}
The part involving dense matrix blocks may be done via the usual block sparse multiplication separately from the low rank part. The latter is done by first applying the row basis tree $\mathcal{V}^T$ to $x$ and involves an upsweep in the $\mathcal{V}$ tree, with multiplication by the explicit basis at the leaves and the transfer matrices up the tree.
This is followed by the application of $\mathcal{S}$ to this product
for all blocks at all levels. Finally, the application of
$\mathcal{U}$ on this intermediate result is done via a downsweep
through the tree using the transfer matrices and the explicit bases at
the leaves. This is the key idea of fast multipole methods and a
similar operations count shows that the complexity is $O(kn)$. In
addition, there is much concurrency in all stages and overlap between
the stages can be exploited for efficient execution on GPUs \cite{boukaram18b}.

Recompression is a basic operation on hierarchical matrices that we rely on for efficiently building and updating them. Consider for example the task of adding a low rank update to $A_{\mathcal{H}^2}$,
\begin{equation}
	\tilde{A}_{\mathcal{H}^2} = A_{\mathcal{H}^2} + X Y^T
\end{equation}
where $X$ and $Y$ are global low rank matrices of size $n \times k'$. This update will evidently affect all blocks of $A_{\mathcal{H}^2}$ at all levels in the hierarchy, since the triple product form of a given block  in the updated matrix may be expressed as:
\begin{equation}
\tilde{A}^l_{ts} =
	\begin{bmatrix} U_t & X_t \end{bmatrix}
	\begin{bmatrix} S_{ts} & 0 \\ 0 & I \end{bmatrix}
	\begin{bmatrix} V_s^T \\ Y_s^T \end{bmatrix}
\end{equation}
where $X_t$ and $Y_s$ are restrictions of $X$ and $Y$ to the index sets $t$ and $s$, 
respectively. The local ranks of the
matrix blocks have now increased to $k + k'$. The leaves of
the basis tree $\mathcal{U}$ and its transfer matrices have also
increased by $k'$. Recompression is the problem of finding new bases
$\mathcal{\bar{U}}$ and $\mathcal{\bar{V}}$ with rank $\bar{k} < k +
k'$ and projecting the matrix $\tilde{A}$ on these bases, to obtain
the compressed representation $\bar{A}$ such that the error introduced
is below the target threshold to which $A$ itself is represented,
$||\tilde{A} - \bar{A}|| \le \epsilon$. The new
bases may be expressed as $U_t R_t$ where $R_t$ is the result of
computing a $QR$ factorization of the block row defined by $t$, and it can be
computed quite efficiently via a downsweep pass through the basis
trees performing $QR$ factorizations only on small dense blocks. $U_t
R_t$ is then truncated to the target accuracy via an SVD or randomized
SVD via an upsweep pass through the basis trees performing SVD
operations only on small dense blocks. All recompression computations
may be done in log-linear time complexity, and there is substantial
parallelism in all of the tree sweep operations, which can be exploited
in manycore architectures. For further algorithmic and implementation details on the
recompression procedure, see \cite{boukaram18b}.

Matrix norms are also required by various operations, and they can be efficiently obtained via fast matrix-vector products or by accessing the low rank blocks. The exact Frobenius norm of a hierarchical matrix with an orthogonal basis can be efficiently evaluated in $O(kn)$ by simply taking the sum of the squares of the Frobenius norms of the coupling matrices of all blocks, $\norm{A}_F^2 = \sum \norm{A_{ts}}_F^2 = \sum \norm{U_t S_{ts} V_s^T}_F^2 = \sum \norm{S_{ts}}_F^2$. On the other hand, $p$-norms can be approximated by a sequence of matrix-vector products using for example an iterative $p$-norm power method \cite{higham2002accuracy}. While the $1$-norm and $\infty$-norm only require a very small number of iterations to converge, the 2-norm can require a relatively large number of iterations. In this case, we terminate the iterations when the norm value has converged to only 2 significant digits. In practice this is a good enough threshold for the relative norm estimation used in this work. 


\subsection{Construction from Hessian-vector products}
\label{sec:hara}
Here we describe how to construct the hierarchical
Hessians {\em ab initio}. In \cite{yu17, yu18}, a method for
constructing a hierarchical representation is presented that relies on
having $O(1)$ access to individual matrix entries. In
\cite{martinsson11}, construction of an HSS approximation is presented
that also requires $O(1)$ access to entries. Unfortunately, in the
case of the Hessians originating from PDE-constrained problems, access to individual
entries in constant time is not generally possible, as these operators are
only available in the form of matrix-vector products.
For example, given an arbitrary vector in $\mathbb{R}^n$, the Hessian
matrix-vector product for inverse problems governed by
(time-dependent) PDEs is computed by solving the (forward-in-time)
forward and (backward-in-time) adjoint equations as well as second order
incremental state-like and adjoint-like PDEs, and by assembling the
resulting product from spatio-(temporal) integrals involving the
state, adjoint, incremental state, and incremental adjoint solution
variables \cite{HinzePinnauUlbrichEtAl09}. Examples of this procedure
for inverse problems governed by several different classes of PDEs are
described in \cref{sec:applications}.

\begin{figure}
	\includegraphics[width=0.3\textwidth]{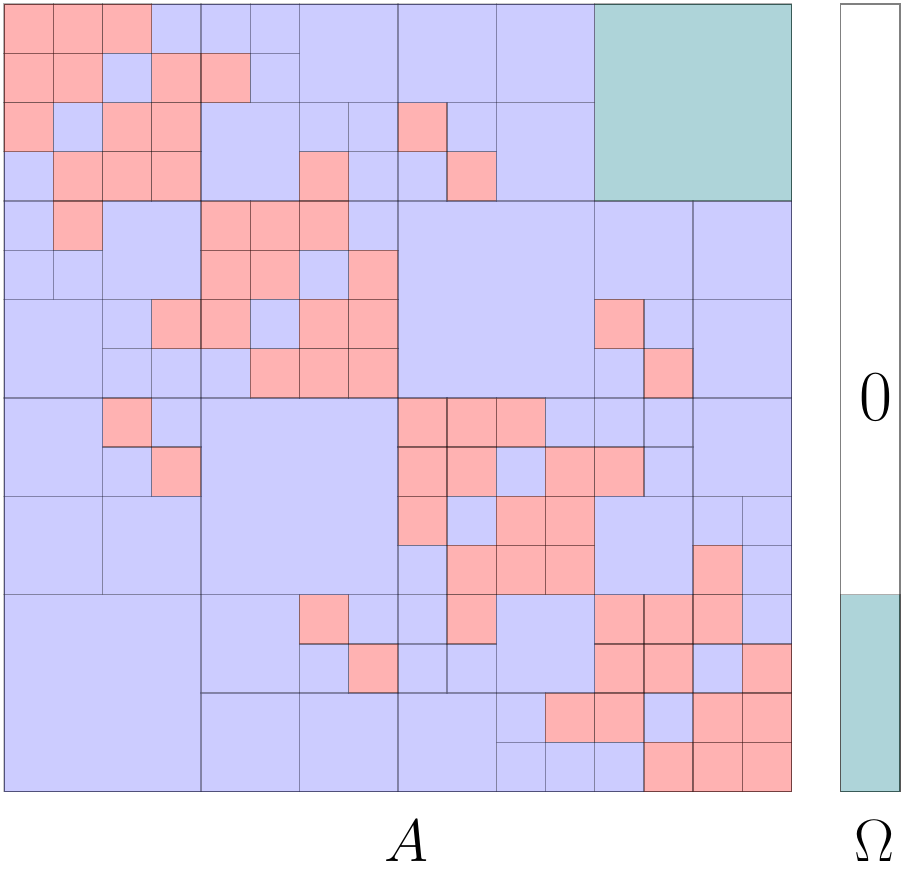}
	\hspace*{.03\textwidth}
	\includegraphics[width=0.3\textwidth]{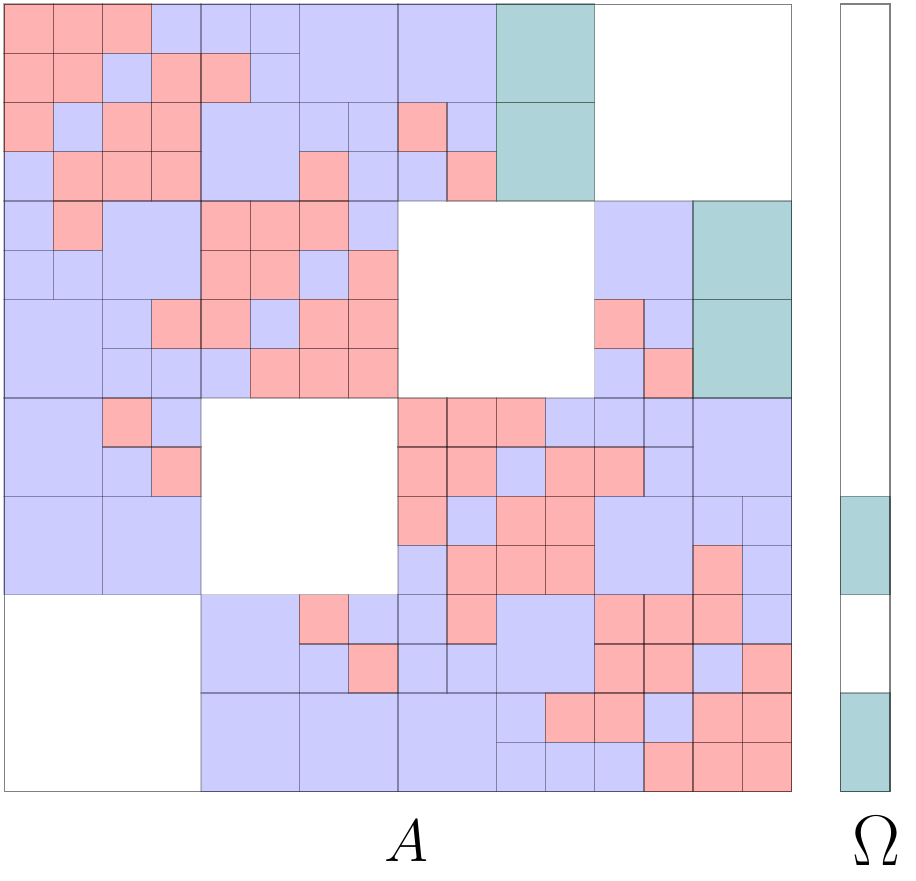}
	\hspace*{.03\textwidth}
	\includegraphics[width=0.3\textwidth]{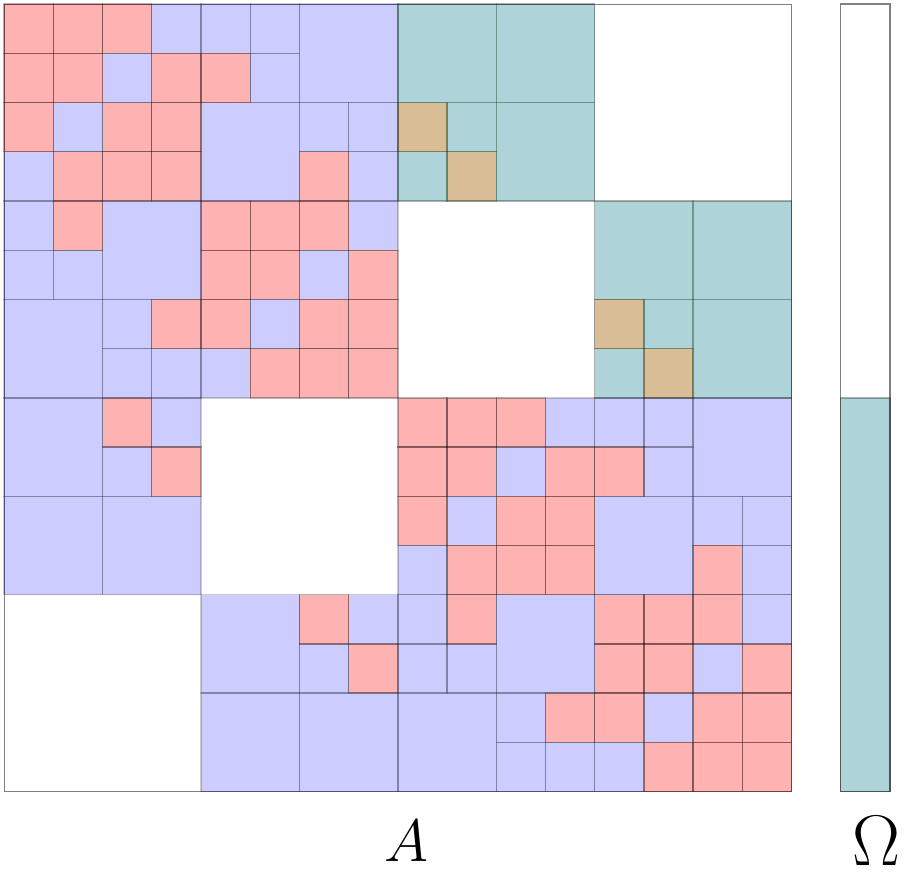}
	\caption{Block sampling patterns. Sampled blocks need not correspond to the block structure of the hierarchical Hessian being constructed as illustrated in the figure on the right.}
	\label{fig:sampling}
\end{figure}

In order to generalize the successful randomized methods of dense
matrix factorizations \cite{halko11} to hierarchical matrices, a
number of procedures \cite{ying11, martinsson16, boukaram19} have been
proposed. As with dense matrices, these procedures sample the matrix
via products with random vectors. The generalization to
hierarchical matrices proceeds by levels, from coarsest to finest,
and uses patterned random vectors chosen to sample particular blocks
of the matrix. \cref{fig:sampling} illustrates how this local sampling works; a column
space for the  low rank representation of the top right block (in the left
figure) is generated by matrix-vector products with random columns $\Omega$
having nonzero values corresponding to the block being sampled.
The resulting product has in its top block row, $(A\Omega)_1$, the desired
column space, and the rest of the product may be discarded. A $QR$
factorization of $(A\Omega)_1$ gives the orthogonal columns basis $U$
and the product $V = A^T [0 \ U^T]^T$ gives the row basis for the block
and its representation $UV^T$. This process is repeated until the approximation 
of the sampled block satisfies an approximation threshold \cite{boukaram19}. 
Similar low rank representation for the
remaining blocks at the same level are performed, and they are then
added as \emph{local} low rank updates to the hierarchical matrix being
constructed. 

The blocks of the next level in the hierarchy are processed in a
similar fashion but by sampling a matrix from which the higher level
blocks previously sampled have been removed, as illustrated in the
middle panel of \cref{fig:sampling}. Note that a single sampling
matrix can be used to generate low rank representations of multiple
blocks simultaneously as long as other blocks do not
interfere with the sampling. Once the low rank block representations
are generated, they are added as local low rank updates and the matrix
is recompressed.  We note here that we may sample larger blocks than
the destination matrix has, as shown in the right panel of \cref{fig:sampling} 
and our current implementation uses a weak admissibility sampling structure.
Low rank representations of these sampling blocks
are generated similarly and then added to
multiple blocks of the destination matrix via local low rank
updates. If the ranks of these blocks are too large there may be an
advantage in sampling with larger blocks as the matrix is peeled in
fewer steps, albeit with every step requiring more matrix-vector products.

For every level $l$ in the matrix, we therefore need $C_l k$ Hessian
vector products, where $C_l$ is a constant that depends on how many blocks can be sampled simultaneously at level $l$,
and $k$ is a representative local rank.
The total number of samples is then $C k \log n$. 
The constant $C$ is small if we choose
large sampling blocks  (which will likely have larger local ranks
$k$), but it is larger if we choose sampling strategies with smaller
blocks. Therefore  tuning the granularity of the sampling, which
affects $C$ and $k$, may be needed, and it depends on the structure of the
Hessian and on the relative costs of the matrix-vector products. It is also worth noting here that on modern hardware architectures, doing $C_l k$ matrix-vector products,
i.e., solving multiple state and adjoint PDEs simultaneously, may not be
much more expensive that solving for a single column, as the computations
involved are  usually memory bound  and a multiple right-hand side approach
will substantially increase the arithmetic intensity and the resulting throughput of the algorithm. As a
result, operations count may not be the best complexity metric to
evaluate performance.

The above procedure may be combined with global
low rank approximations and may be accelerated if a Hessian
approximation is available, e.g., for a different parameter from a
previous iteration. For example, in situations where the data available for inversion
is limited, the Hessian of the data misfit has a small finite
dimensional range space, and the global Hessian rank may be small
relative to its size. In these cases, an initial global small-rank
approximation $BB^T$ of the misfit Hessian may be computed by the
usual dense randomized sampling methods. If the global rank
approximation does not reach the desired target accuracy, the
hierarchical construction can be done on the residual $(H - B B^T)$
matrix, and $B B^T$ may then be added as a low rank update to it. This
will likely require a smaller total number of matrix-vector products for the sampling
relative to directly building the hierarchical Hessian.

A similar procedure can be used if a Hessian
approximation is available, e.g., from a previous optimization
iteration for a different parameter. Denoting the Hessian
approximation by $\tilde H$, one can construct a hiearchical
approximation (or a global low rank approximation) for the difference
$H-\tilde H$. Alternatively, provided $\tilde H$ is invertible and we
have an inverse square root  $\tilde H^{-1/2}$ available, one could
compute an approximation to $(\tilde H^{-1/2} H \tilde H^{-1/2}-I)$.
In the same spirit, regularization matrices, which are generally
sparse but have full rank, may be readily added via local low rank updates
after the hierarchical misfit Hessian is constructed to obtain an explicit full
Hessian representation that can then be efficiciently operated on. We discuss iterative 
inversion methods, which may also be adapted to square root and inverse square 
root computations, in \cref{sec:inverse}. Clearly, other linear algebraic methods
for these computations are also possible once the compressed matrix representation 
of the Hessian is available. 

\section{Illustrative Applications}
\label{sec:applications}

In this section, we discuss the effectiveness of the hierarchical
representation of Hessians for several inverse problems, governed by a
time-dependent diffusion equation in one spatial dimension
(\cref{sec:diffusion}), and in two space dimensions by a steady-state
advection-diffusion equation (\cref{sec:adv-diff}), a
frequency-domain wave equation (\cref{sec:wave}), and a
time-dependent low-frequency Maxwell equation
(\cref{sec:em}). The last two examples are prototypes for
geophysical inverse problems with seismic and controlled-source
electromagnetic (CSEM) modalities, respectively.

The problems we consider are infinite-dimensional inverse problems,
i.e., we infer parameter fields. We use a deterministic inverse
problem approach, i.e., we use an optimization formulation and employ
regularization to cope with the inherent ill-posedness. Note that
under certain assumptions, the resulting PDE-constrained optimization
problem can also be interpreted in a Bayesian context, in which case
the minimizer corresponds to the maximum a posteriori parameter
estimate.  Upon discretization with finite elements (details on the
discretization is given for each problem individually), these problems
results in large-scale optimization problems, which we solve using an
inexact Newton-conjugate gradient descent method with linesearch \cite{EW96,NocedalWright06}. This
requires availability of gradients and of Hessian-vector products,
which are computed using adjoint methods as summarized in each example
separately below. Similar methods are common in PDE-constrained
optimization, \cite{HinzePinnauUlbrichEtAl09, BorziSchulz12}, and have
successfully been used in large-scale optimization formulations of
inverse problems, e.g., \cite{Bui-ThanhGhattasMartinEtAl13,
  EpanomeritakisAkccelikGhattasEtAl08, ZhuLiFomelEtAl16,
  IsaacPetraStadlerEtAl15, HesseStadler14}. Here, we mainly focus on
the approximations of the Hessian, which are always evaluated at the
parameter field found as the solution of the optimization problem. In
particular, we study the cost in terms of the number of Hessian
applications (in PDE-constrained optimization this is usually the
dominating cost) required for the construction of H-matrix
approximations, and compare it with the cost of global low rank
approaches. We also study the influence of the number of observations
and properties of the governing PDEs on the compressibility of Hessian
matrices.



\subsection{Density inversion in 1D time-dependent diffusion equation} 
\label{sec:diffusion}

\begin{figure}[!ht]
\begin{tabular}[c]{cc}
	\begin{subfigure}[c]{0.6\textwidth}
		\hspace*{4pt}
		\includegraphics[width=.85\textwidth]{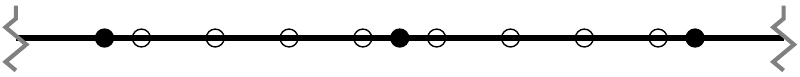} \\
		\includegraphics[width=.85\textwidth]{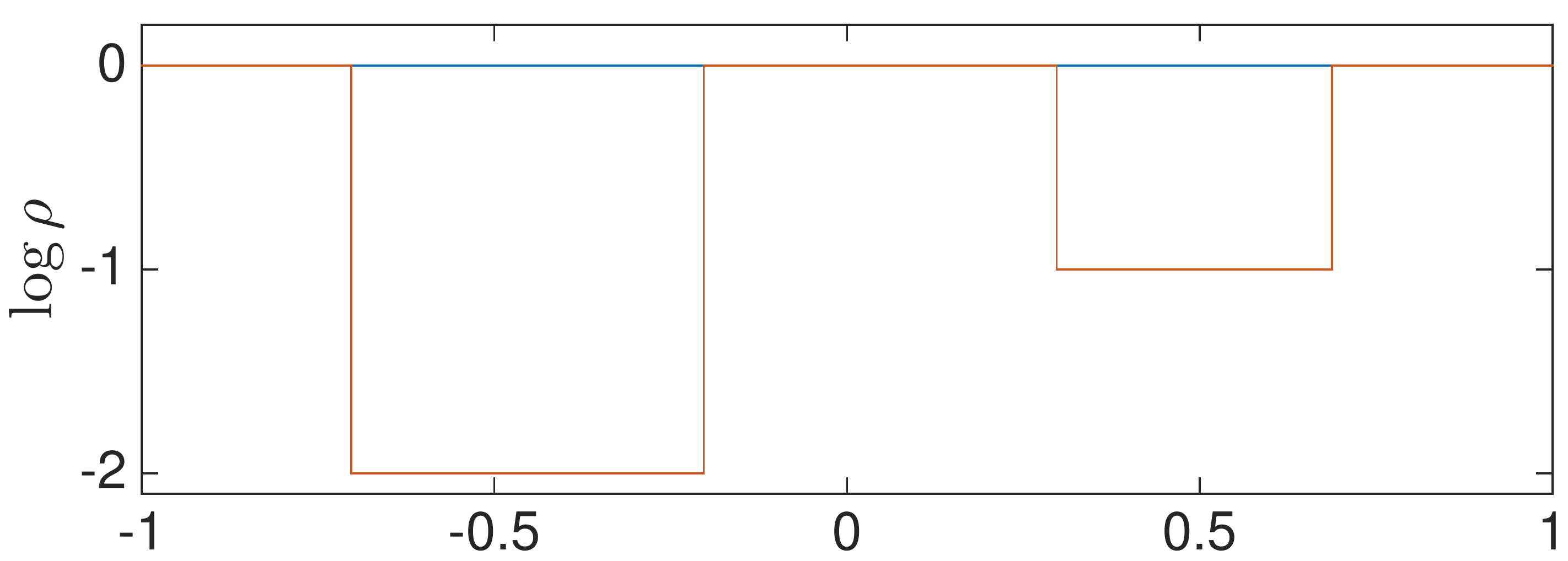}
	\end{subfigure} 
	&
	\begin{subfigure}[c]{0.3\textwidth}
		\hspace*{6pt}
    	\includegraphics[width=.75\textwidth]{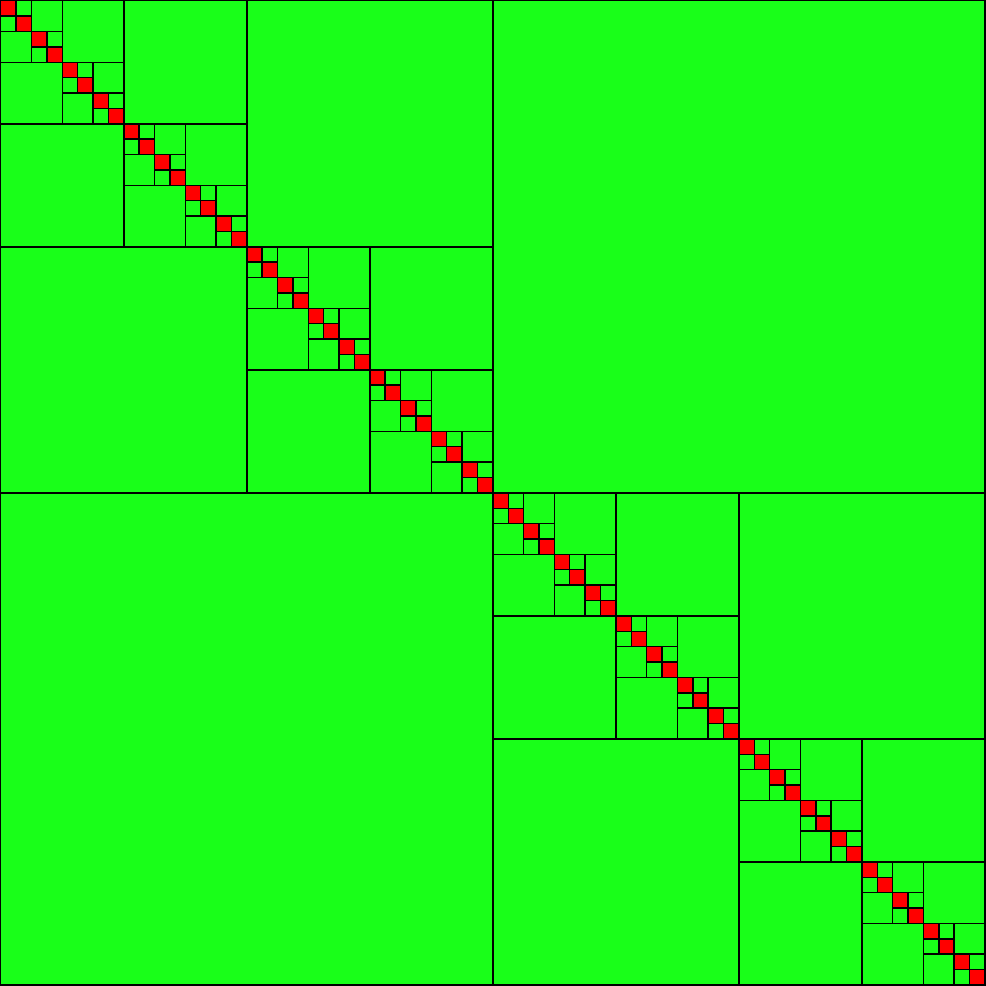} 
	\end{subfigure} 	
\end{tabular}
\caption{(left) Sources ($\bullet$), receivers ($\circ$), and target profile; (right) Hessian structure.}	
\label{fig:setup1d}
\end{figure}

We consider a one-dimensional domain where we seek to invert
for a spatially-varying porosity coefficient $\rho(x)$ in
the governing diffusion equation $\rho(x) \partial_t u - \partial_x ^2
u = w(x, t)$. We use 3 sources and 8 receivers placed as shown in
\cref{fig:setup1d}. The point sources produce a Ricker wavelet
time history input $W$ and the signals are recorded at the receivers
until they have effectively dissipated at a final time $T>0$.

The medium to be recovered is in the interval $-1 \le x \le 1$ and has
sharp discontinuities (see \cref{fig:setup1d}), prompting the use of a
(smoothed) total-variation regularization such that the optimization
formulation is:
\begin{equation}
	\underset{\rho:[-1,1]\to\mathbb R}{\mbox{minimize}}\:\:
        J(\rho)
        := \tfrac{1}{2} \sum_s \sum_r \int_0^T (u_s(x_r, t) - d_r(t))^2 dt 
	+ \alpha \!\int_{-1}^1\! \sqrt{|d_x \rho|^2 + \beta} \, dx, 
\end{equation}
  {where for all sources $s$, $u_s(x,t)$ is the
    solution to}
\begin{equation}\label{eq:diffusion-state}
  \begin{aligned}
        \rho (x) \, \partial_t u_s - \partial_x ^2 u_s &= \delta(x - x_s) W(t - t_0), \quad
	                    &&|x| < \infty, \ 0 \le t \le T, \vspace*{0.25em} \\
	 u_s(\pm \infty, t) &= 0, \hspace*{2.65cm}  &&\text{{\small (Dirichlet\ BC)}}  \vspace*{0.25em}\\
	 u_s(x, 0) &= 0. \hspace*{4.5cm} &&\text{{\small (IC)}}  \\
  \end{aligned}
\end{equation}
Here, the first term of the objective is the least squares time misfit
summed over all sources and receivers. The second term is
a total-variation regularization functional that allows sharp
discontinuities to be recovered, where $\alpha>0$ is the
regularization weight and a small $\beta>0$ ensures that the objective
is differentiable.

For the computations below, $\rho(x)$ is discretized in $[-1,1]$ using
$n$ linear elements. A finite element discretization of the governing
equation with the same mesh is used in $[-1,1]$ and extended outside
with constant $\rho = 1$ to apply the homogeneous
Dirichlet boundary conditions at a sufficiently far distance. An
implicit second-oder time integration scheme is used for solving the
semi-discrete equations.

The computation of the gradient may be done at the cost of two PDE
solutions, a forward and an adjoint. The continuous
form of the gradient of $J$, its functional Fr\'echet derivative, is
\begin{equation}
	\frac{\delta J}{\delta \rho} (x) = \sum_s \int_0^T \partial_t u_s(x, t) p_s(x, t) dt
	\label{eq:grad1d}
\end{equation}
where, for each source $s$, $u_s$ is the solution of the state equation
\eqref{eq:diffusion-state}, and $p_s$ the solution of the
corresponding adjoing equation that must be solved backwards-in time:
\[
\begin{array}{cc}
%
	\text{(adjoint)} & 
	\begin{array}{c}
		-\rho(x) \, \partial_t p_s - \partial_x ^2 p_s = - \sum_r \delta(x - x_r) (u_s(x_r,t)-d_r(t)) \vspace*{1pt} \\
		q(\pm \infty, t)=0; \quad q(x, T) = 0.
	\end{array}

\end{array}
\]

To compute a hierarchical representation of the Hessian, we
need the ability to compute products of the Hessian of the misfit term
with vectors generated as described in \cref{sec:Hrep}. This may be done at the cost of solving, for each
source $s$, two additional
forward-like and adjoint-like PDEs. The continuous form of the
Hessian-vector product, i.e., the product of the Fr\'echet Hessian
with a model perturbation $\nu(x)$ is: 

\begin{equation}
	\frac{\delta^2 J}{\delta \rho^2} \nu =  \sum_s \left( \int_0^T \partial_t u(x, t) q(x, t) dt + 
         \int_0^T \partial_t v(x, t) p(x, t) dt \right) + \frac{\delta^2 R}{\delta \rho^2} \nu,
\label{eq:Hv1d}
\end{equation}
where the last term denotes the Hessian of the TV regularization $R$,
and $v_s$ and $q_s$ are the solutions of the following incremental or
second order forward and adjoint equations:
\[
\begin{array}{lc}
	\text{\small (2\textsuperscript{nd} order forward)} & 
	\begin{array}{c}
		\rho(x) \, \partial_t v_s - \partial_x ^2 v_s = - \nu(x) \, \partial_t u_s(x, t), \vspace*{.1pt} \\
		v_s(\pm \infty, t) = 0; \quad v_s(x, 0) = 0,
	\end{array}
	\vspace*{6pt}\\\

	\text{\small (2\textsuperscript{nd} order adjoint)} & 
	\begin{array}{c}
		\hspace*{-8pt}
		-\rho(x) \, \partial_t q_s - \partial_x ^2 q_s = - \sum_r \delta(x - x_r) v_s(x_r,t) - \nu(x) \, \partial_t p_s(x, t), \vspace*{.1pt} \\
		q_s(\pm \infty, t)=0; \quad q_s(x, T) = 0.
	\end{array}

\end{array}
\]

\begin{figure}[!ht]
	\begin{center}
		\includegraphics[width=0.31\textwidth]{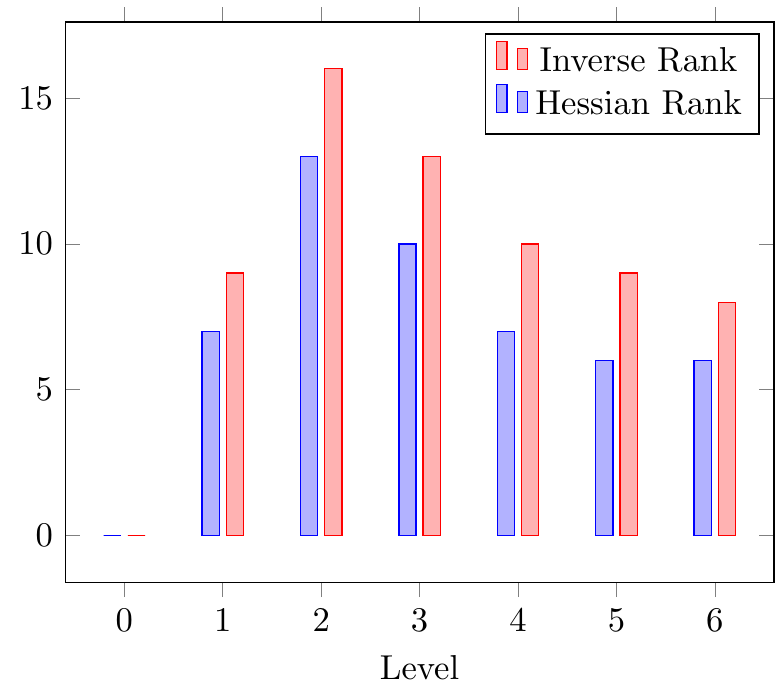}		
		\includegraphics[width=0.34\textwidth]{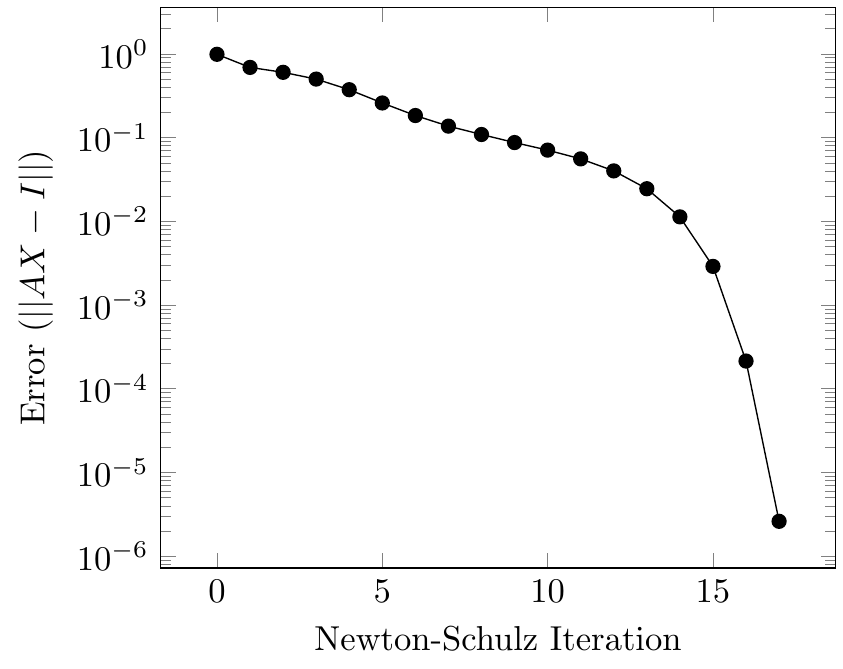}
		\includegraphics[width=0.33\textwidth]{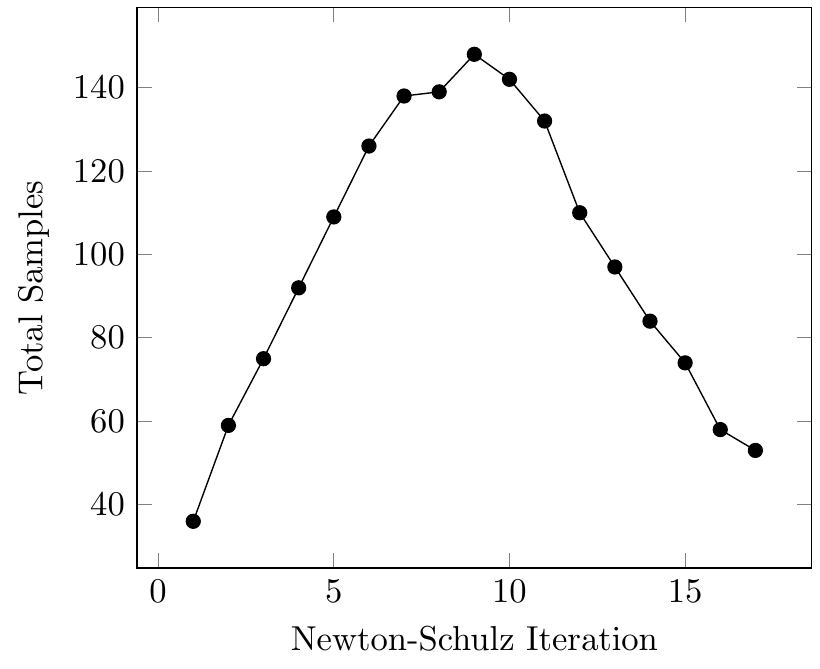}
	\end{center}
	\caption{(left) Local ranks of Hessian and its inverse. (middle) Convergence history of NS. (right) Number of samples needed to construct the matrix over the iterations.}
	\label{fig:ns1d}
\end{figure}

The structure of the $\mathcal{H}$-representation of the Hessian is depicted in
the right panel of \cref{fig:setup1d}. We used a weak
admissibility partitioning of the hierarchical matrix as appropriate
for a one-dimensional spatial domain. There are $2^l$ symmetric blocks
pairs at level $l$, each of size $n/2^l$. For this example of size
$n=2048$, we use 6 levels in the hierarchy, stopping the matrix
refinement at blocks of size $m = 32$. We construct the Hessian to a
relative accuracy of about $10^{-6}$ in the 2-norm. The local ranks of
the resulting matrix are shown in the left panel of
\cref{fig:ns1d}. The ranks plotted are the maximum local ranks for all
off-diagonal blocks for every level. Even though the Hessian is full
rank because of the regularization term, we note that its local ranks
are small relative compared to its size even when a relatively tight
tolerance of $10^{-6}$ is used in the construction.

Given this explicit representation of the Hessian, we can invert it to
produce an explicit inverse.
The convergence behavior of Newton Schulz starting from a scaled identity
initial iterate to a relative accuracy of $\epsilon = 10^{-6}$ is
shown in \cref{fig:ns1d}. The number of samples needed to construct
the hierarchical approximation of the iterates is shown as well. In
this simple 1D setting, we used a static threshold and a constant
accuracy in the iterates. The increase in the required samples, and
the corresponding local ranks, of the intermediate iterates is
noteworthy. The left panel of \cref{fig:ns1d} shows the local ranks of
the inverse along with the Hessian local ranks for comparison. There
is little change in the local ranks in the inversion operation.

\subsection{Source inversion in stationary advection-diffusion} 
\label{sec:adv-diff}
We consider a linear source inversion problem, i.e., we infer the right
hand side source
$m(x)$ in an advection-diffusion-reaction equation from point
observations of its solution on a domain $\Omega\subset \mathbb
R^2$. This amounts to the following optimization problem:
\begin{equation}
	\underset{m:\Omega\to\mathbb R}{\mbox{minimize}} \:\:
                J(m) := \frac{1}{2\sigma^2} \sum_{r} (u(x_r) -
                   d_r)^2+ {R}(m), \\
\end{equation}
\noindent where $u$ is the solution of
\begin{equation}\label{eq:adv-diff-state}
\begin{aligned}
	 -{\rm div}(\kappa \nabla u) + {\bf v}\cdot \nabla u + cu &=
         m && \text{in } \Omega,\\
	 u &= 0  && \text{on } \Gamma_D,  \vspace*{0.25em}\\
	 \kappa \frac{\partial u}{\partial n} &= 0  && \text{on } \Gamma_N. \\
\end{aligned}
\end{equation}
Here, the measurements $d_r$ are assumed to contain additive
independent and identically distributed Gaussian noise with variance
$\sigma$, and ${R}(m)$ is a regularization term.  The
coefficients $\kappa$, ${\bf v}$ and $c$ represent the diffusivity, the
advective velocity and the reaction constant, respectively, and $\Gamma_D,
\Gamma_N \subset \partial \Omega$ is a splitting of the
boundary $\partial\Omega$, where we impose Dirichlet or Neumann
boundary conditions, respectively.
\label{sec:source}
\begin{figure}[!ht]
	 \centering
	 \includegraphics[height=0.22\textwidth]{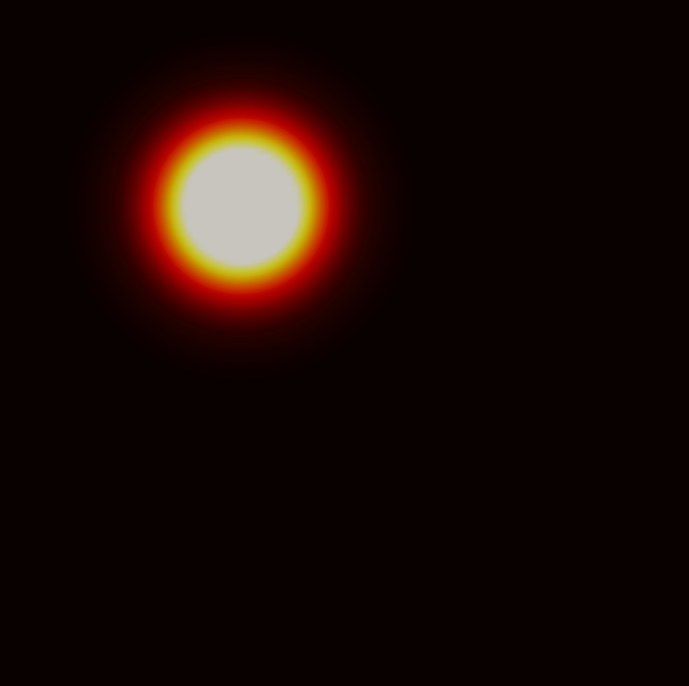}
         \includegraphics[width=0.03\textwidth]{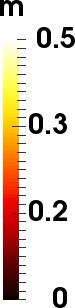}
         \hfill
	\includegraphics[height=0.22\textwidth]{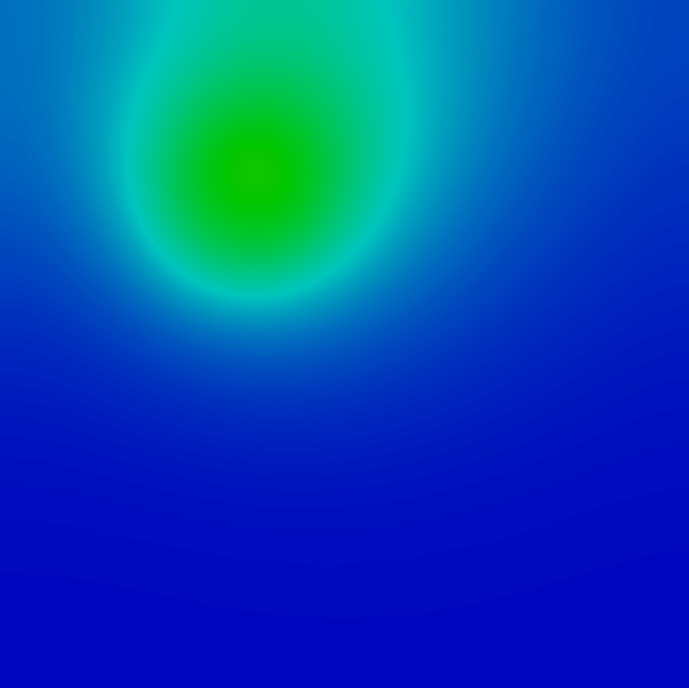} 
	\includegraphics[height=0.22\textwidth]{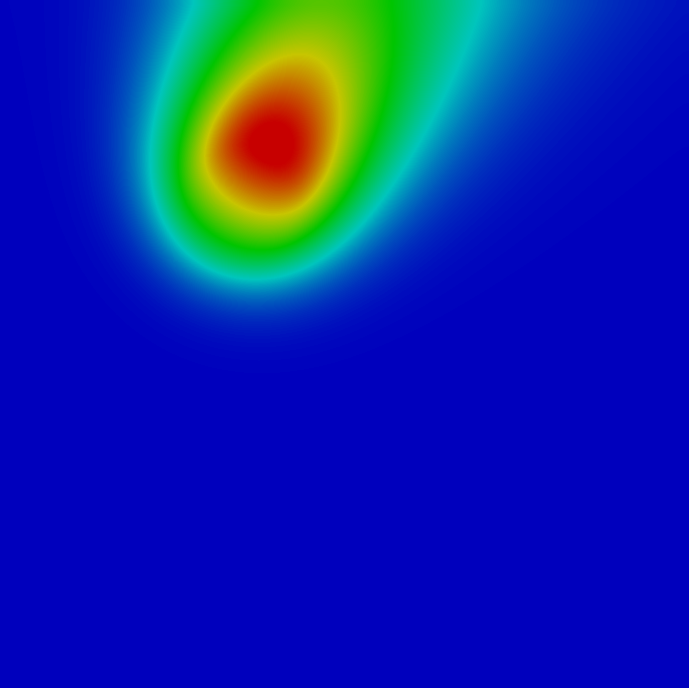} 
	\includegraphics[height=0.22\textwidth]{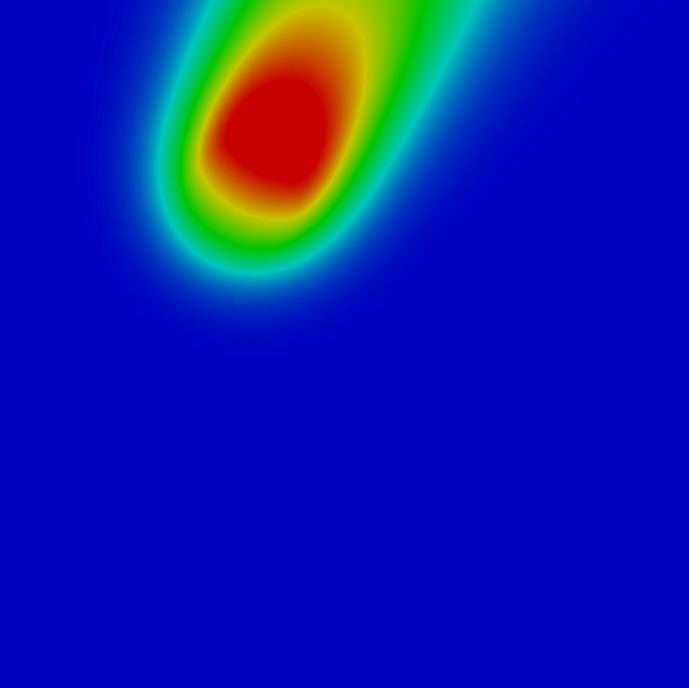}
        \includegraphics[width=0.03\textwidth]{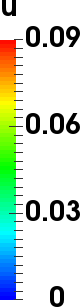}\\[1ex]
	\caption{Dependence of the solution for different diffusion
          parameters $\kappa$: Shown on the left is the source $m$
          entering on the right hand side of
          \eqref{eq:adv-diff-state}. The right three plots show the
          solutions of the state equation for $\kappa=10^{-1},
          10^{-2}, 10^{-3}$ (from left to right).  Depending on
          $\kappa$, the problem is diffusion- or advection-dominated.}
	\label{fig:source_map}
\end{figure}

For our study of the Hessian, we use $\Omega=[0,1]^2$, the
physical parameters $c=0.5$ and ${\bf v} = (x_1,x_2)^T$, and a noise
level of $1\%$. To
study the influence of the diffusion and the number of observations,
we use different values of $\kappa \in
\{10^{-3},10^{-2},10^{-1}\}$ and different numbers of observation
points, $N_{obs} \in \{ 250, \dots,
16000\}$.

We discretize $\Omega$ using triangles based on a uniform grid of size
$128 \times 128$, and use linear continuous finite elements to
discretize the parameter $m$. The presence of strong advection in the
problem mandates the use of stabilization for the state and adjoint
variables. We employ the first-order SIPG scheme \cite{ArnoldBrezziCockburnEtAl02,
  riviere2008discontinuous} with upwinding
\cite{ayuso2009discontinuous} for the forward and adjoint problems. Our implementation is based on
FEniCS \cite{LoggMardalEtAl2012a} and hIPPYlib, an open-source library
providing scalable adjoint-based algorithms for PDE-based inverse
problems \cite{VillaPetraGhattas2018}. The true source and the state solutions corresponding to the different values of the diffusion coefficient are shown in  \cref{fig:source_map}. 


\begin{figure}[!ht]
	\includegraphics[width=0.32\textwidth]{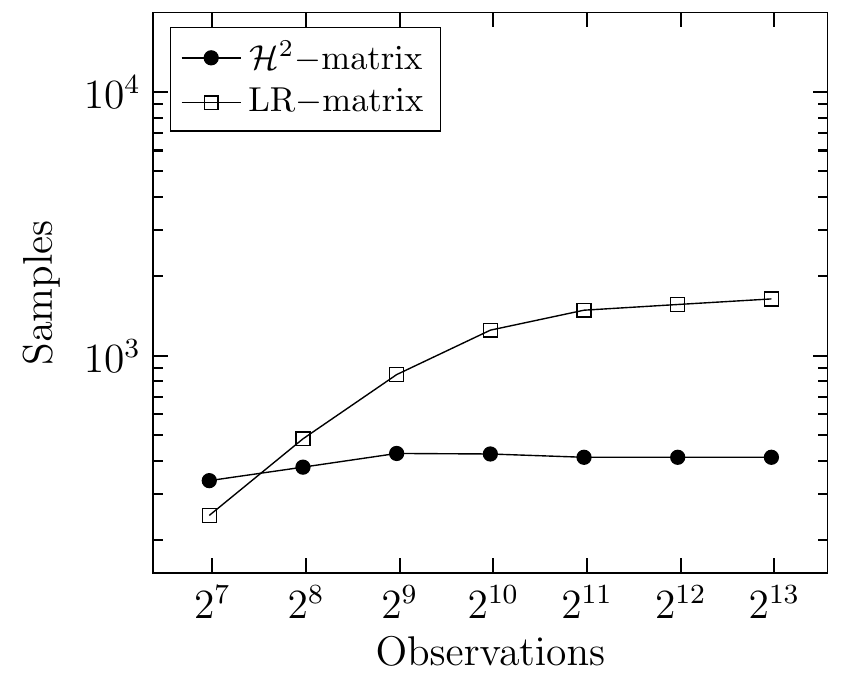}	
	\includegraphics[width=0.32\textwidth]{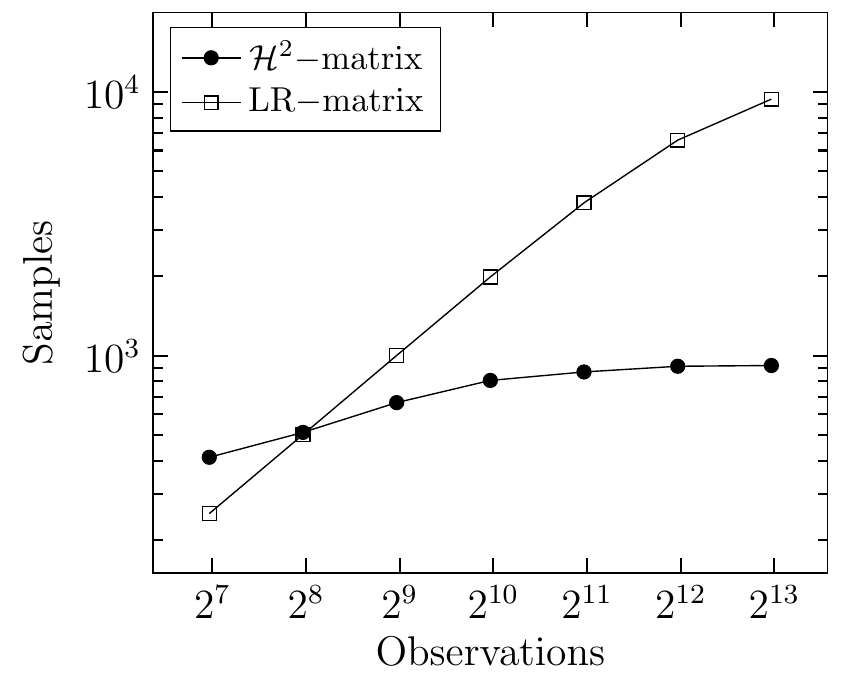}	
	\includegraphics[width=0.32\textwidth]{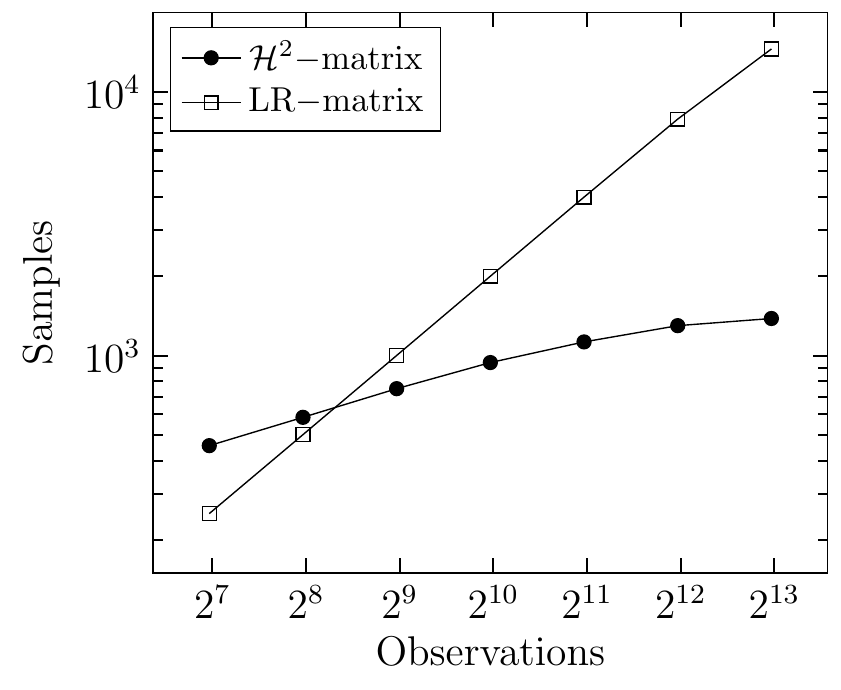}
	\caption{The $y$-axis shows the number of Hessian-vectors
          products (samples) needed for constructing hierarchically
          and globally low rank approximations of the misfit Hessian
          to an accuracy of $\epsilon = 10^{-4}$ for increasing number
          of observations plotted on the $x$-axis. Shows are
          results for $\kappa=10^{-1}$ (left), $10^{-2}$ (middle), and
          $10^{-3}$ (right).}
	\label{fig:source_samples}
\end{figure}

In our study of the compressibility of the Hessian and the associated
computational cost, we only consider the misfit Hessian, i.e., neglect
the Hessian of the regularization term. This allows us to compare $\mathcal{H}$-matrix
approaches with the more commonly used global low rank approach, which
would suffer from inclusion of the Hessian of the regularization,
which typically has full rank.  Note that since this is a linear inverse
problem, this misfit Hessian does not depend on the parameter $m$.
\cref{fig:source_samples} shows the number of Hessian-vector products
needed for constructing the misfit Hessian to a relative spectral
accuracy of $10^{-4}$ for three different values of the diffusion
coefficient $\kappa$ and an increasing number of observations.
%
For small numbers of observations, the misfit Hessian has a fast
decaying spectrum and a globally low rank representation can be
achieved at a lower cost than with the hierarchical approach, measured
in the number of Hessian-vector products. However, as more
observations are incorporated, the Hessian applications required for
the global low rank approximation grows rapidly.  In contrast, the
cost for the hierarchical compression only increases mildly with the
number of observations as it is insensitive to the global rank.
We also observe that larger diffusion coefficients $\kappa$ limit the
increase of computational cost for increasing number of observations
for both, the hierarchical and the global low rank
approximations. This is a consequence of the fact that diffusion
limits the amount of fine-scale information that can be recovered in
the inversion, and thus more observations do not have a significant
effect on the rank of the misfit Hessian.

\begin{figure}[!ht]
\centering
\begin{subfigure}[b]{0.33\textwidth}
	\raisebox{11pt}{\includegraphics[width=\textwidth]{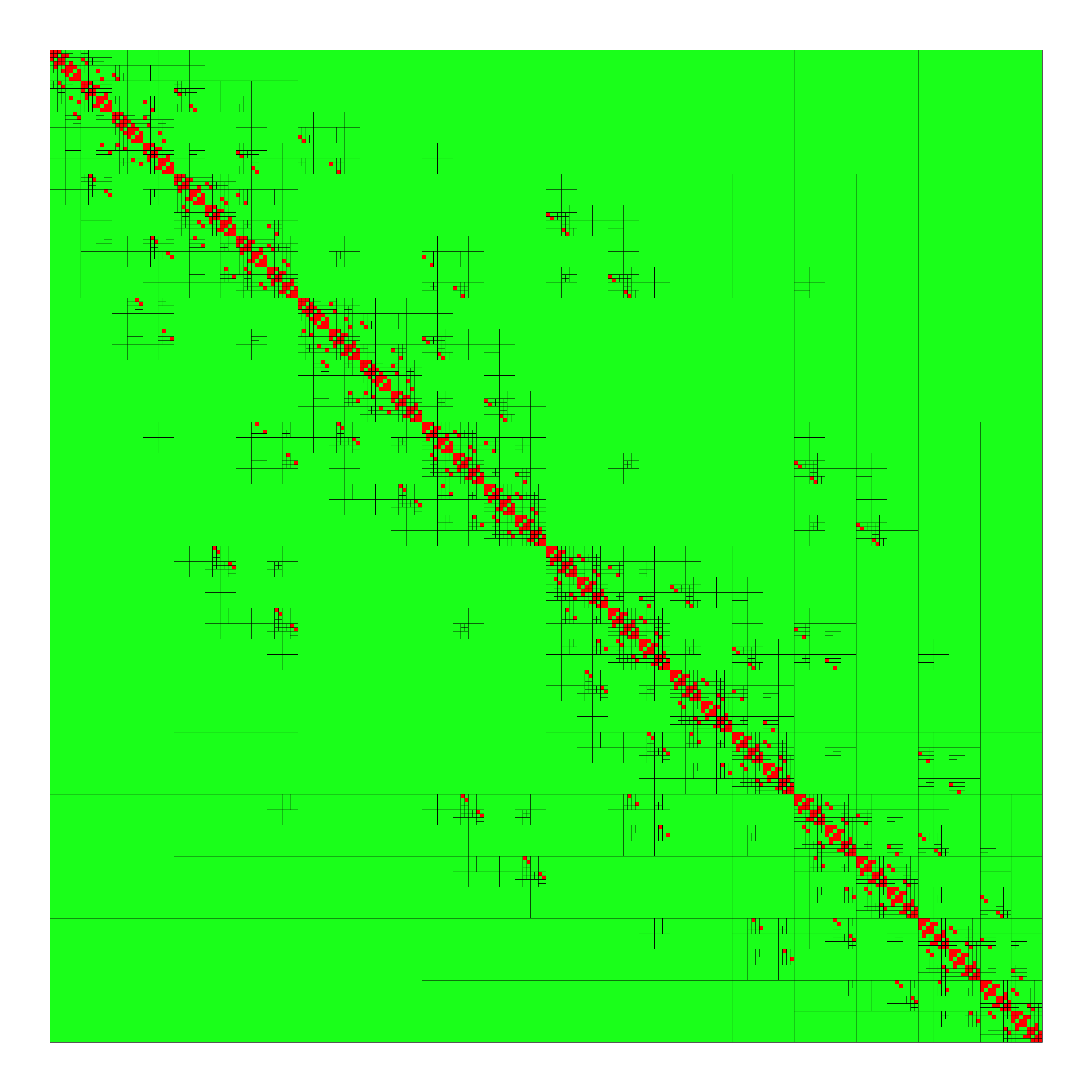}}
\end{subfigure}
\hspace*{20pt}
\begin{subfigure}[b]{0.49\textwidth}
	\includegraphics[width=\linewidth]{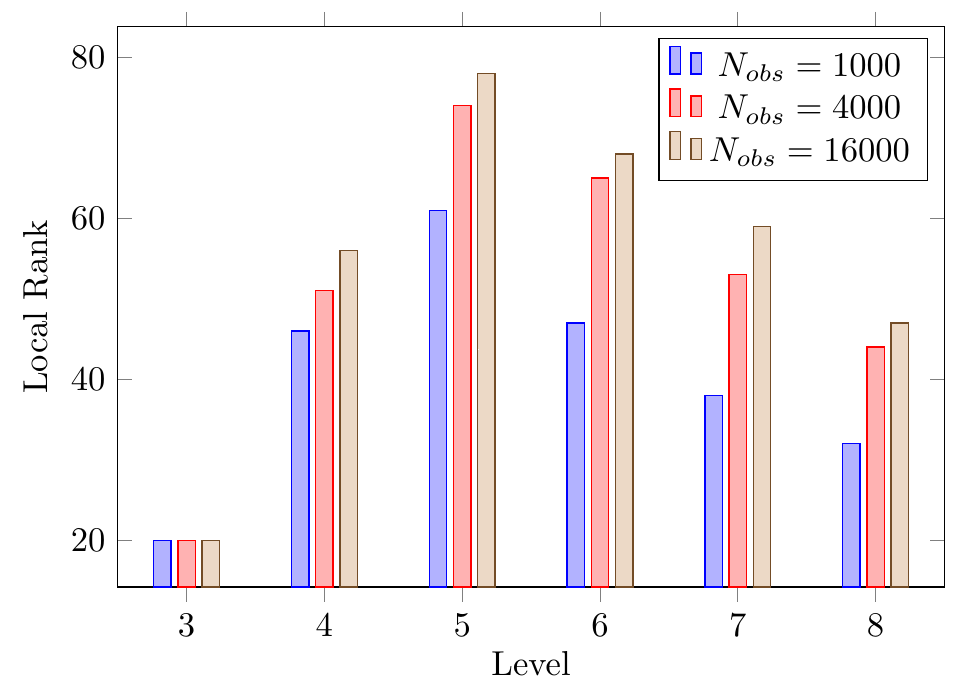}
\end{subfigure}
\caption{8-level structure of Hessian (left) and corresponding maximum local
  ranks per level for different number of observations and
  $\kappa=10^{-3}$ (right).}
\label{fig:source_ranks}
\end{figure}

The relative insensitivity of the number of samples to the data
dimension is reflected in the fact that the blocks of
the resulting $\mathcal{H}^2$ have small local ranks that do not grow
much even as the global rank of the matrix
increases. The left panel of \cref{fig:source_ranks} shows the structure of the
constructed matrix which uses a node ordering of the mesh obtained from a KD-tree binary space partitioning method, with leaf clusters of size 64, and an admissibility condition that resulted in the refined matrix partitioning shown.
The right panel of \cref{fig:source_ranks} shows how the block ranks grow much slower than the global Hessian rank as the data dimension of the problem grows. 


\subsection{Frequency-domain wave equation inversion}
\label{sec:wave}

%
%

We consider the inverse acoustic wave propagation problem in frequency
domain governed by the Helmholtz equation in a bounded domain $\Omega\subset
\mathbb R^2:$
\begin{equation}
	\Delta u + \omega^2 q u = f,  \quad \text{ in } \Omega,
	\label{eq:sec52-forward-problem}
\end{equation}
where $u$ is the time-harmonic pressure, $f$ is the source and
$w$ is the angular frequency. One corresponding inverse problem is to
infer for the squared slowness,  $q(x) = 1/c^2(x)$, where $c(x)$ is
the speed of sound. We assume that the observations are given as the
frequency-weighted normal derivative of the time-harmonic pressure at
the receivers, located at all discretization points $x_r$ on the top boundary of the
domain.

Assuming homogeneous Dirichlet boundary conditions and introducing the
linear observation operator $B$ defined by $B( u)(x_r) = \omega^2
\nabla u(x_r) \cdot n$, the deterministic inverse problem with an
appropriate regularization term $R(q)$ can written as follows:
\begin{equation}
\begin{aligned}
	&\underset{q:\Omega\to\mathbb R}{\text{minimize } } \:\: {J(q) := \tfrac{1}{2} \sum_r  ( B( u)(x_r) - d_r)^2
	                                          + R(q)}, \\
	&\quad \text{where $u$ solves }  \text{\eqref{eq:sec52-forward-problem} with
        } u = 0 \text{ on } \partial \Omega.
\end{aligned}
\end{equation}
In order to be able to compare with global low rank methods, we do not
include the regularization term in the gradient and
Hessian expressions below, nor include them in our compression
experiments.
The continuous gradient of the misfit objective $F$ at a parameter
$q_0$ is given by:
\[
\frac{\partial F}{\partial q} (q_0) =  \omega^2 q_0 u_0  v,
\]
where $u_0$ is the solution of the state problem \eqref{eq:sec52-forward-problem} with slowness $q_0$ and $v$ solves the adjoint problem:
\[
\begin{array}{l}
	{\displaystyle\Delta v + \omega^2 q_0 v = -\sum_r  (B( u)(x_r) - d_r)  \, \text{ in } \Omega, \quad v = 0 \, \text{ on } \partial \Omega.}
\end{array}
\]
In our numerical experiments, we use the Gauss-Newton approximation
$H^{\text{GN}}$ of the Hessian, which neglects terms involving the
adjoint variable.  This approximation is commonly used in practice since,
differently from the Hessian, it is guaranteed to be positive
semidefinite. In matrix form, it is given by
\[ H^{\text{GN}} = C^TA^{-T}B^TBA^{-1}C,
\]
where $A$ stems for the discretization of the state operator, and
$C$ from the discrete representation of the partial derivative of the
PDE with respect to the parameter.

\begin{figure}
\centering
\includegraphics[width=0.75\textwidth]{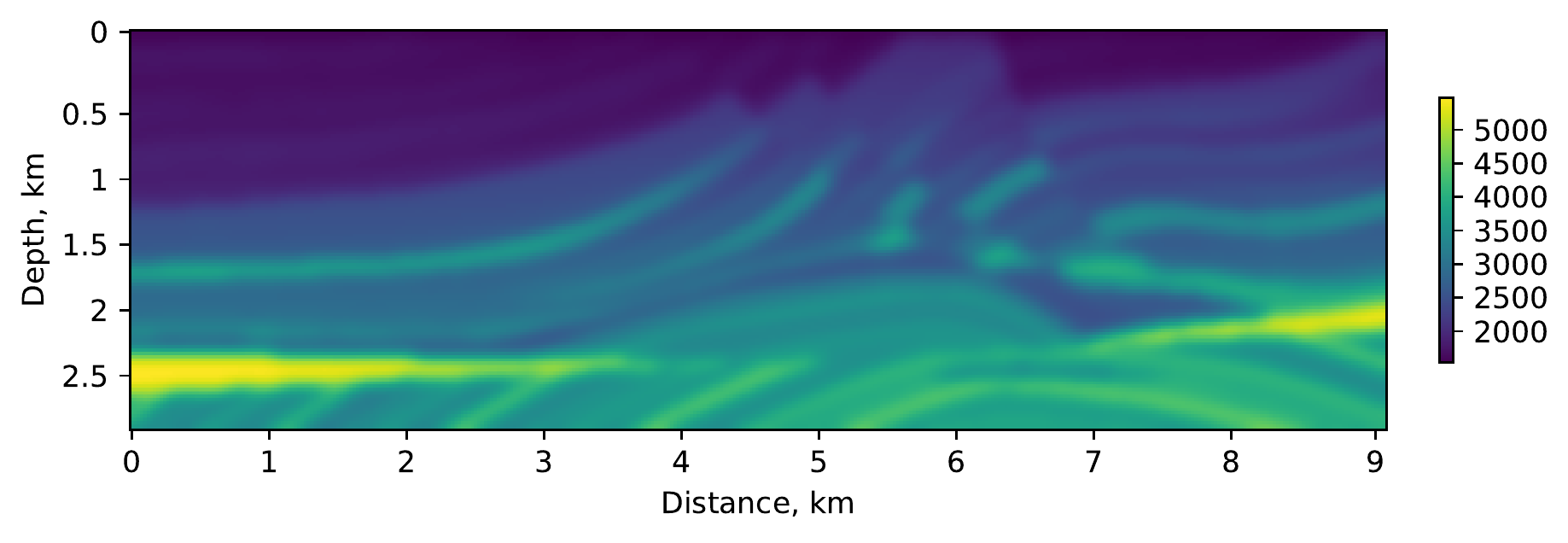}
\caption{Smoothed Marmousi model used to generate synthetic
  observations.}
\label{fig:wave_truth}
\end{figure}
Our numerical test is set up using the 2D acoustic Marmousi model
\cite{versteeg1994marmousi}, a benchmark model for seismic
inversion. The initial ``hard'' model is smoothed using a Gaussian
filter to produce the result shown in \cref{fig:wave_truth}, which we
use to generate synthetic observations.
We consider a triangulation of the
the rectangular domain $\Omega = [0, 9192]\times [0,
  2904]$ m$^2$ in order to obtain about ten grid points per wave
length.
The state, adjoint and parameter variables are discretized using continuous
piecewise linear Lagrangian finite elements.  Since we use the same finite element basis for
the state $u$ and the parameter $q$, the $ij$ entry of C reads as
$$C(u_0,q_0)_{ij} = \omega^2 \int_\Omega \phi_j u_0 \phi_i ~dx.$$
The right-hand side term
in \eqref{eq:sec52-forward-problem} is given as a sum of sources,
located at a depth of 10m, spaced evenly every 12.5m.  We further
incorporate Dirichlet boundary condition along the top boundary of the
domain and PML absorbing boundary conditions
\cite{berenger1994perfectly} on the other three boundaries. As in the
previous problem, our implementation uses
FEniCS and hiPPYlib
\cite{LoggMardalEtAl2012a, VillaPetraGhattas2018}.


\begin{figure}[!ht]
	\centering
	\includegraphics[width=0.65\textwidth]{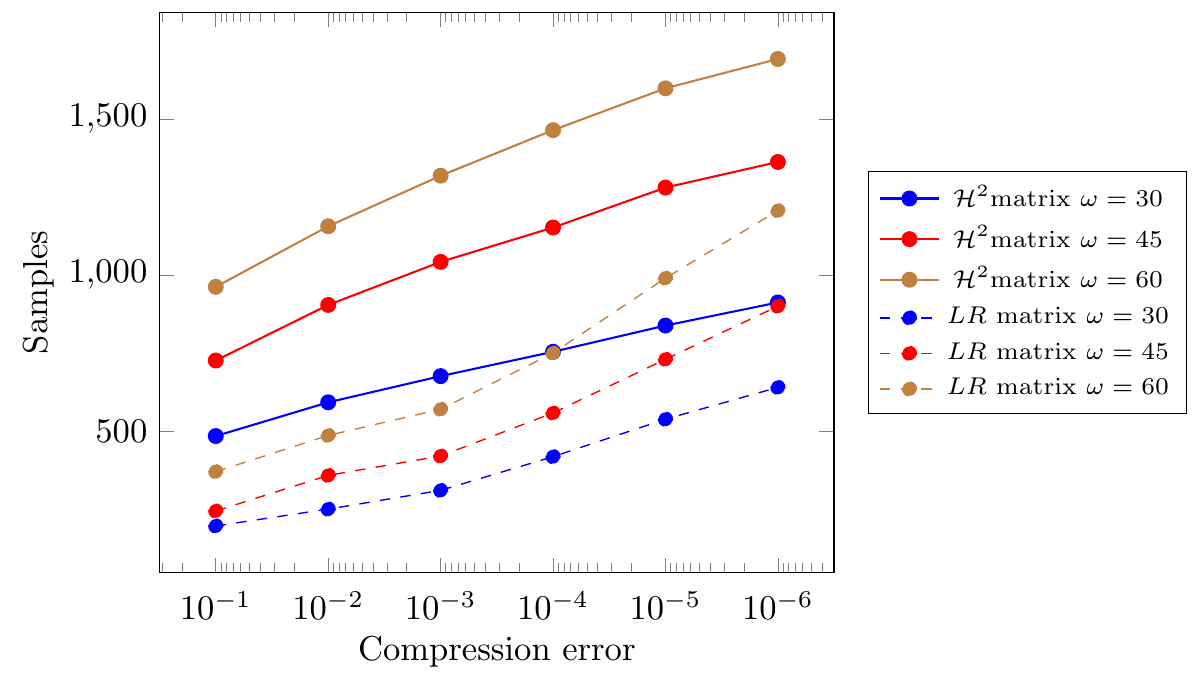}
	\caption{Number of Hessian applications (samples) needed to
          construct the Gauss-Newton Hessians for
          different frequencies and different accuracies for
          hierarchical and global low rank approximations.}
	\label{fig:wavesamples}
\end{figure}


\begin{figure}[!ht]
	\centering
	\begin{tabular}[c]{cc}
		\includegraphics[width=0.3\textwidth,valign=t]{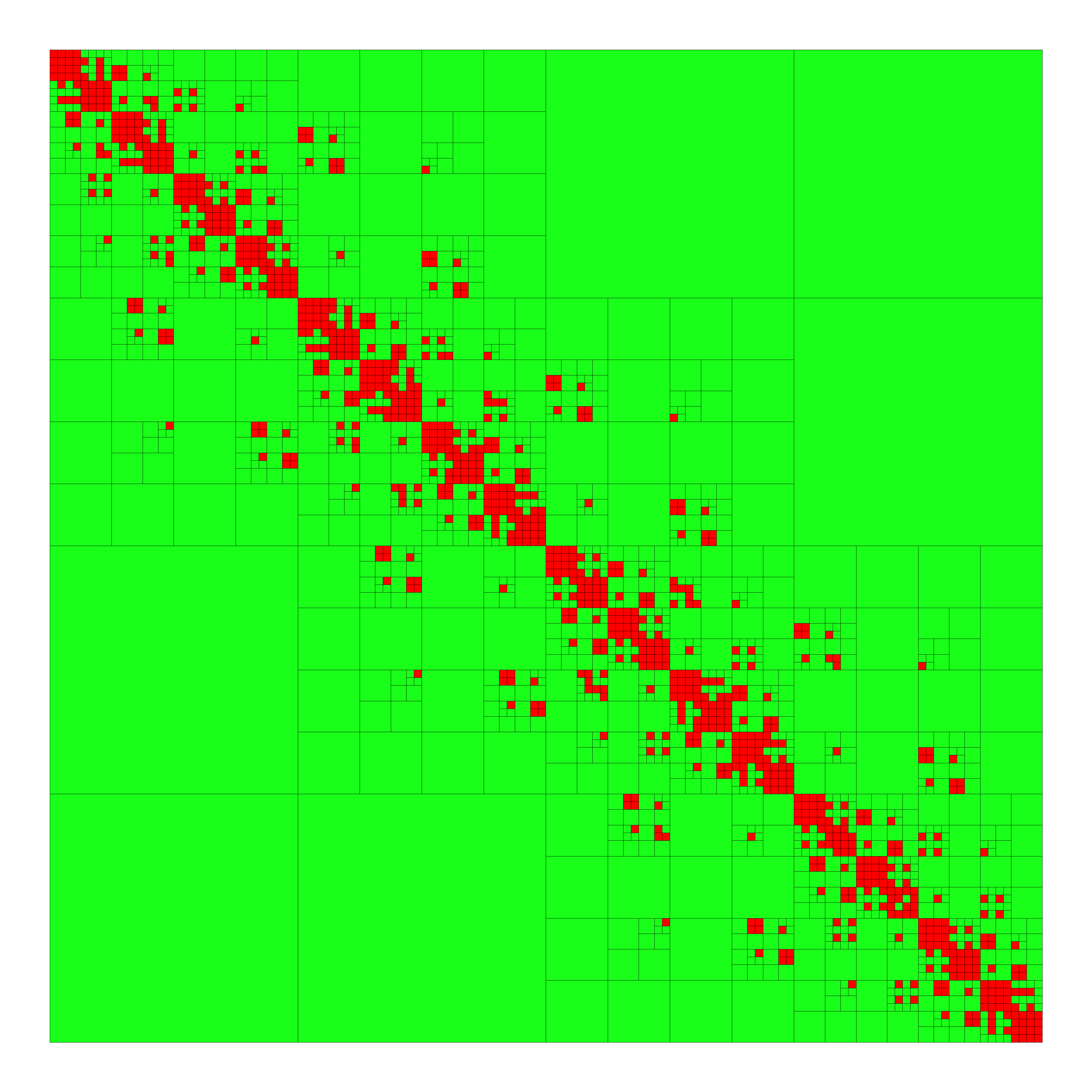} &
		\includegraphics[width=0.55\textwidth,valign=t]{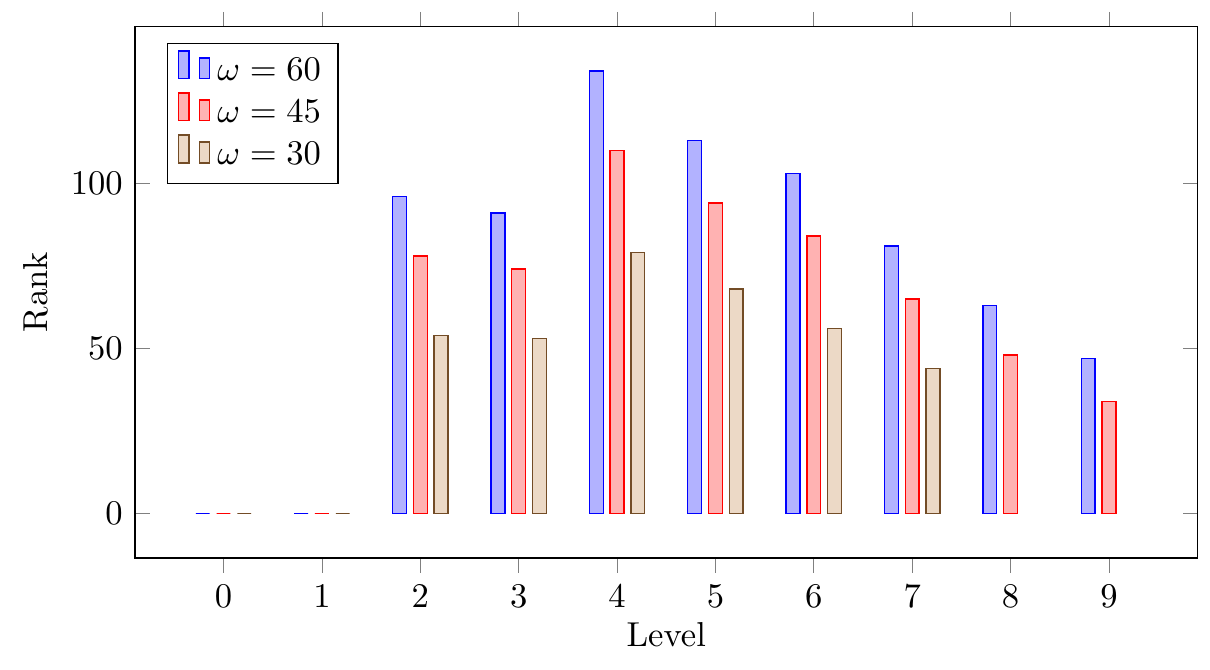}
	\end{tabular}
	\caption{Structure and max local ranks per level of Hessian at
          $\epsilon = 10^{-4}$.}
	\label{fig:waveranks}
\end{figure}

\cref{fig:wavesamples} shows the number of samples needed to
construct the Gauss-Newton Hessian evaluated at the converged
parameter values for a range of frequencies and to accuracies in the
range $10^{-6} \le \epsilon \le 10^{-1}$ measured in the relative
2-norm of the matrix. The sizes of the Hessians at the three
frequencies $w = 30$, $45$, and $60$ are $n = 8{,}211$, $18{,}164$,
and $33{,}269$, respectively, growing with the square of the frequency
to maintain about ten grid points per wave length. For a desired
accuracy, the number of samples needed for the construction of the
Hessian shows a slightly sublinear growth with frequency, and for a
given frequency, a slow growth in the number of required samples with
$|\log \epsilon|$. Due to the fairly limited information, coming from
a single supersource input, the global low rank of the misfit Hessian
requires fewer Hessian-vector products than the hierarchical
representation. However as the accuracy increases the global grows
rapidly, particularly at higher frequencies. We expect that with more
data incorporated in the inversion (using multiple independent
sources, for example) the total number of required samples will grow
slower for the hierarchical compression than for the global low rank
approximation.

\cref{fig:waveranks} shows the structures and the local ranks of the
resulting compressed Hessians. A leaf size of 64 was used in the
construction with an admissibility condition that produces the
structure shown. The matrix structure for the larger problems
corresponding to $\omega=45, 60$ is visually similar to the matrix
structure of the smaller problem $\omega=30$ but is refined by two
more levels.  The resulting local ranks (maximum block rank) per level
are shown for a relative accuracy of $\epsilon = 10^{-4}$. We note the
relatively mild sublinear growth in the local ranks with frequency.


\subsection{Transient controlled-source electromagnetic inversion}
\label{sec:em}
In this section we consider a problem arising in transient
controlled-source electromagnetism, see, e.g.,
\cite{Simoncini2013,Haber14}. The deterministic optimization approach
to infer the spatially varying electrical conductivity field
$\vb*{\sigma}$ from observation data minimizes the
least squares misfit between these observations and the response
predicted by the model, i.e.,
\begin{equation}
	\underset{\vb*{\sigma}:\Omega\to\mathbb R^{d\times
            d}}{\mbox{minimize }} \: {\displaystyle J(\vb*{\sigma}) :=
          \tfrac{1}{2} \sum_r \sum_s \int_0^T\!
          (\delta(\vb{x}_r)\vb{E}_s - \vb*{d}_{rs})^2 dt + \alpha R(\vb*{\sigma}) = F(\vb*{\sigma}) + \alpha R(\vb*{\sigma}).}
\label{eq:maxwell_misfit}
\end{equation}
Here, $\vb{E}_s$ is the electric vector field which is assumed to satisfy
Maxwell's equations in the low frequency regime
\begin{equation}
\begin{aligned}
\vb*{\sigma} \, \partial_t \vb{E}_s + \nabla \times (\mu^{-1} \nabla \times \vb{E}_s) &= - \vb{J}_s &&\text{ on $\Omega \times (0,T]$}, \\
\vb{E}_s \times \vb{n} &= 0 &&\text{ on $\partial\Omega \times [0,T]$},\\
\vb{E}_s &= 0 && \text{ on $\Omega \times \{0\}$},
\label{eq:maxwell_fwd}
\end{aligned}
\end{equation}
where $\mu>0$ is the constant corresponding to the magnetic permeability of free space,
$\vb*{d}_{rs}$ is the
measured time-varying response at a receiver located at $\vb{x}_r$ corresponding to
a given source term $\vb{J}_s$. In our study, we consider a horizontal
electric dipole point source at $\vb{x}_s$ oriented along the $x$-axis
direction, which leads to $\vb{J}_s = \delta(\vb{x}_s) \vb{e_1}
\partial_t W(t)$, where $W$ is a Ricker wavelet of the type
\cite{Carcione2010}
\[
W(t) := \left( a-\frac{1}{2}\right)e^{-a}, \quad a = \left(\frac{\pi(t-t_s)}{t_p}\right)^2, \quad t_s = 1.4 ~ t_p,
\]
and $\delta(\cdot)$ is the Dirac delta function.
The continuous gradient of $F$ (we do not
include the regularization in the expressions below) at a parameter
$\vb*{\sigma}_0$ may be computed as:
\[
\frac{\partial F}{\partial \vb*{\sigma}} (\vb*{\sigma}_0) = \sum_s \int_0^T \partial_t \vb{E}_{s} \cdot  \vb{P}_{s} ~dt,
\]
where $\vb{E}_{s}$ is given by solving the state equation
\eqref{eq:maxwell_fwd} for the $s$-th source and with conductivity
$\vb*{\sigma_0}$, and $\vb{P}_{s}$ is obtained by solving the adjoint PDE
\begin{equation}
\begin{aligned}
-\vb*{\sigma}_0 \, \partial_t \vb{P}_{s} + \nabla \times (\mu^{-1} \nabla \times \vb{P}_{s}) &= - \sum_r (\delta(\vb{x}_r) \vb{E_{s}}(t) - \vb*{d}_{rs}(t)) && \!\!\text{on $\Omega \times [0,T)$},\\[-1ex]
\vb{P}_{s} \times \vb{n} &= 0 &&\!\!\text{on $\partial\Omega \times [0,T]$},\\
\vb{P}_{s} &= 0 && \!\!\text{on $\Omega \times \{T\}$}.
\label{eq:maxwell_adjoint}
\end{aligned}
\end{equation}
Taking another variation, the product of the Hessian of $F$
at a parameter $\vb*{\sigma}_0$, with a model perturbation
$\vb*{\nu}:\Omega\to\mathbb R^{d\times d}$ is computed as
\begin{equation}
\frac{\partial^2 F}{\partial \vb*{\sigma}^2}(\vb*{\sigma}_0) ~ \vb*{\nu} =  \sum_s \int_0^T \partial_t \vb{E}_{s} \vb{Q}_{s} ~dt + \sum_s \int_0^T \partial_t \vb{F}_{s} \vb{P}_{s} ~ dt,
\label{eq:maxwell_H}
\end{equation}
where $\vb{E}_{s}$ and $\vb{P}_{s}$ are given by solving
\eqref{eq:maxwell_fwd} and \eqref{eq:maxwell_adjoint}, respectively,
and $\vb{F}_{s}$ and $\vb{Q}_{s}$ are obtained from the solution of
the incremental forward problem
\begin{equation}
\begin{aligned}
\vb*{\sigma}_0 \, \partial_t \vb{F}_{s} + \nabla \times (\mu^{-1} \nabla \times \vb{F}_{s}) &= - \vb*{\nu}\partial_t \vb{E}_{s} && \text{ on $\Omega \times (0,T]$},\\
\vb{F}_{s} \times \vb{n} &= 0 &&\text{ on $\partial\Omega \times [0,T]$},\\
\vb{F}_{s} &= 0 && \text{ on $\Omega \times \{0\}$},
\label{eq:maxwell_tlm}
\end{aligned}
\end{equation}
and the incremental adjoint problem
\begin{equation}
\begin{aligned}
-\vb*{\sigma}_0 \, \partial_t \vb{Q}_{s} + \nabla \times (\mu^{-1}
\nabla \times \vb{Q}_{s}) &= - \sum_r \delta(\vb{x}_r) \vb{F}_{s} -
\vb*{\nu} \partial_t \vb{P}_{s} && \text{ on $\Omega \times [0,T)$},\\[-1ex]
\vb{Q}_{s} \times \vb{n} &= 0 &&\text{ on $\partial\Omega \times [0,T]$},\\
\vb{Q}_{s} &= 0 && \text{ on $\Omega \times \{T\}$}.
\label{eq:maxwell_soa}
\end{aligned}
\end{equation}
%

The Hessian matrix-vector products are computed by using the software framework PETScOPT \cite{PETScOPT}.
We discretize the state and ajoint variables with hexahedral Nedelec vector finite elements
using the open-source library MFEM \cite{MFEM}, which supports adaptive
mesh refinement and higher-order elements. For the parameter space, we use standard Lagrange elements.
For time integration of \eqref{eq:maxwell_fwd}, \eqref{eq:maxwell_adjoint},
 \eqref{eq:maxwell_tlm}
and \eqref{eq:maxwell_soa}, we use the Crank-Nicholson scheme as implemented in the PETSc library
\cite{PETSc}.
For all the PDE considered, at each time step, we need to solve a poorly
conditioned linear system with the matrices
\[
\frac{1}{\Delta t} M + A, \quad M_{ij} = \int_\Omega \vb*{\phi}_j \cdot \vb*{\phi}_i ~d\vb{x}, \quad A_{ij} =\int_\Omega \mu^{-1} \nabla \times \vb*{\phi}_j \cdot \nabla \times \vb*{\phi}_i ~d\vb{x},
\]
where $\vb*{\phi}_i$ is the $i$-th Nedelec finite element basis
function for the Nedelec
element space and $\Delta t$ the time step. These solves
represents a key challenge for scalability. We use
%
conjugate gradients preconditioned with
the Balancing Domain Decomposition by Constraints solver
\cite{Zampini2016,Zampini2017,Zampini2017b}
combined with a robust choice of the initial guess based on a reduced basis
approximation scheme \cite{Zampini2011}.
For the reported experiments, this strategy leads, remarkably, to an
average of less than one linear iterations per time step.
\begin{figure}
\centering
\includegraphics[width=0.48\textwidth,valign=t]{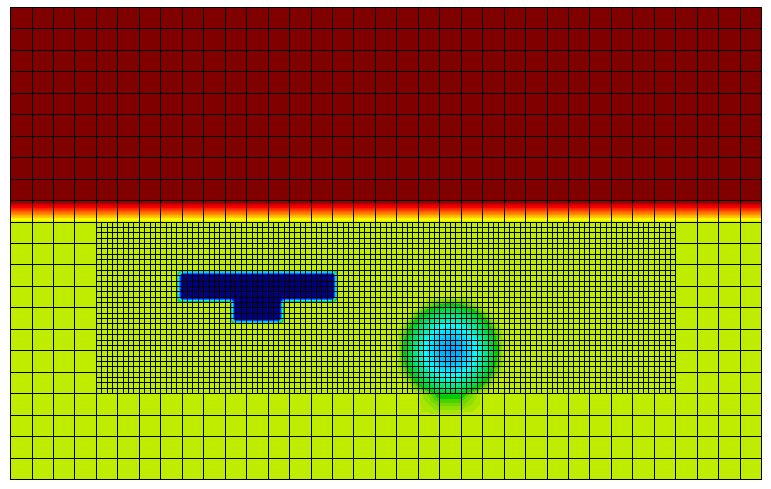}
\includegraphics[width=0.4\textwidth,valign=t]{maxwell_samples}\\[1ex]
\caption{Shown on the left is the truth parameter field to be recoverd
  from CSEM data, together with the locally refined discretization
  overlaid on material properties, which contain two anomalies (in
  blue). The right figure shows the total number of Hessian
  vector samples needed for hierarchical compression as a function of
  accuracy for different regularization parameters $\alpha$ of the TV
  regularizer.}
\label{fig:maxwell_truth}
\end{figure}

To illustrate the hierarchically low rank nature of the Hessians given
in \eqref{eq:maxwell_misfit}, we generate synthetic data for a
two-dimensional isotropic tensor $\vb*{\sigma}$ whose spatial
distribution is given in \cref{fig:maxwell_truth}; the simulated time
is 2 seconds, with a constant time step of 0.01 seconds.  The values
used for the conductivity field, reported in logarithmic scale in the
figure, are 3.0 S/m (water), 0.1 S/m (sediments), 0.01 S/m (salt), and
0.001 S/m for the T-shaped anomaly region.  The non-conforming mesh
used to discretize the conductivity is also shown. The setup uses an
array of 14 equispaced sources located at the water/sediments
interface and an array of 18 receivers placed right in between the
water/sediments interface and the non-conforming mesh interface.
The region of interest for the optimization process
extends to the first layer of elements surrounding the non-conforming
interface, for a total of 2381 optimization degrees of
freedom. For the optimization, we used a primal-dual total variation regularizer
\cite{ChanGolubMulet99}, which preserves sharp gradients in the
coefficient values close to the data recording location, and
guarantees an almost optimal line-search process.

Results for Hessian compressibility reported here consider the full Hessian as
given in \eqref{eq:maxwell_H}.
\cref{fig:maxwell_truth} shows the number of samples needed to
construct the Hessian at convergence of the inversion process. 
We note two characteristics of hierarchical representations in this
numerical experiment. The first is the log dependence of the number of
samples (and local ranks, not shown here) on the desired accuracy. The
second is the relative insensitivity of the local ranks to the amount
regularization. The global rank of the matrix does not have a direct
effect on the local ranks.



\begin{figure}
\includegraphics[width=\textwidth]{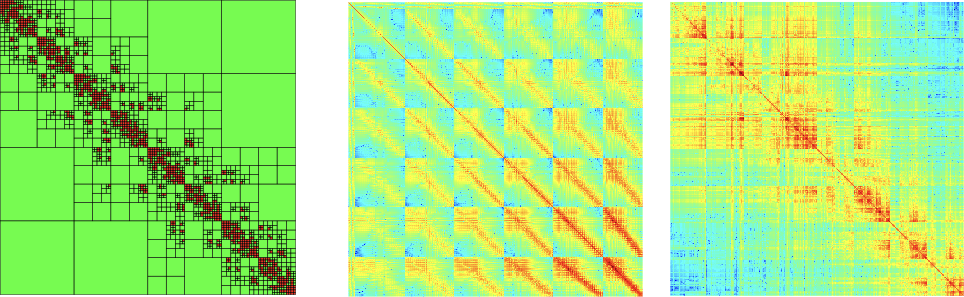}
\caption{
  Shown in the left figure is the structure of the hierarchical
  Hessian using the KD-tree ordering. Middle and right figures
show the Hessian entries with lexicographic and
  KD-tree ordering of mesh nodes, respectively. Color reflect
  magnitude of absolute value, red is largest and blue is smallest.
}
\label{fig:maxwell_h}
\end{figure}

\cref{fig:maxwell_h} shows the structure of the hierarchical low rank
Hessian approximation.
Note the importance of ordering the matrix in a
spatially-aware manner so that indices that are close correspond to
nodes that are in spatial proximity. The figure depicts the
entries of the Hessian when the nodes are ordered lexicographically
and re-ordered using a KD-tree binary spatial partitioning. The
reordering of the optimization degrees of freedom, results in larger
blocks that are ``smooth'', a necessary characteristic for an
efficient hierarchically low rank representation.

%
%
%


\section{Conclusions}
\label{sec:conclusions}

In this paper we presented a hierarchical matrix representation of
Hessians as a viable and practical representation for storing and
manipulating second order information for inverse problems governed by
PDEs, which we expect will be applicable to a broader class of
PDE-constrained optimization problems. We have shown the increasing
superiority of hierarchical matrix approximation of the data misfit
Hessian over a low rank approximation as the data become more
informative for four inverse problems governed by diffusion, transport,
and acoustic and electromagnetic wave propagation. Since the eigenvalues of the data 
misfit Hessian typically decay more slowly in higher dimensional problems, 
we expect the hierarchical approximations to continue to be effective for 
inverse problems in three dimensions, and plan to study their performance in future work.

The primary advantages of the hierarchical representation are that it
provides a tunable accuracy approximation, it has a small and
asymptotically optimal memory footprint, and it may be generated from
Hessian-vector products that are generally available via state and
adjoint PDE solves. Our numerical results have shown that the Hessians
that appear in applications do, as predicted, admit a hierarchical
structure and can therefore be represented with blocks at different
levels of granularity that have small bounded rank. Once the
hierarchical Hessian approximation is computed, it can also be
operated on very efficiently, purely algebraically, for performing
operations such matrix-vector multiplications, local and global low
rank updates, inversion, and square roots. 

The machinery of hierarchical Hessians is general, algebraic, and
effective in storing and manipulating the formally dense
Hessians arising in a variety of PDE-governed inverse problems. Our
initial results on moderately-sized problems are encouraging, and we
intend to study the computational advantages and effectiveness of
hierarchical Hessians at large scale in future work.

\appendix


\section{Iterative Methods for Inverse Hessian Approximation}
\label{sec:inverse}


Iterative methods for computing the matrix inverse have a long history dating back to \cite{schulz33}. Under reasonably general assumptions, they can be shown to be globally convergent methods, and require $2 \log \kappa + \log \log (1/\epsilon)$ number of iterations to converge to an $\epsilon$ accuracy, where $\kappa$ is the condition number of the matrix being inverted \cite{soderstrom74}.

Iterative inversion methods have not been very popular (but see
\cite{sanders15})
because at face value their computational cost is larger than direct
factorization-based inversion methods as they require two $O(n^3)$
matrix multiplications per iteration when used with dense
storage. This calculus changes with hierarchical matrices, as the
product of two such matrices may be obtained in log-linear time
$O(kn\log n)$ \cite{boukaram19}. This fast matrix multiply
makes iterative inversion methods attractive here, even with memory
considerations aside.

\begin{algorithm}[!ht]
\caption{Stabilized Newton-Schulz for inverting $A$ starting with $X_0$}
\label{alg:ns}
\begin{algorithmic}[1]
\Procedure{newton\_schulz}{$A$, $X$, $X_0$, $\epsilon$}
	\State $X = X_0$
	\While{$\norm{AX-I}_2 > \epsilon$}
		\State $X = (2I - X A) X$
		\If{stabilization\_on}  \Comment{\emph{when unwanted $\sigma$ become small} }
			\State $X = X A X$
		\EndIf
	\EndWhile
\EndProcedure
\end{algorithmic}
\end{algorithm}

A basic stabilized Newton-Schulz (NS) iteration
shown in \cref{alg:ns}. For symmetric positive-definite matrices, a scaled identity $X_0 = I /
||A||_\infty$ may be used as the starting iterate and results in
globally convergent iterations. Newton-Schulz converges only linearly
in the early iterations before entering the rapid quadratic
convergence regime. Its convergence can be sped up with improved
initial guesses and convergence accelerations \cite{pan91}. In the
current optimization context, the inverse from a previous iterate
provides a natural starting point and this warm-start reduces the
number of iterations compared to an agnostic starting point,
as we illustrate below. The basic NS iteration is numerically
stable and even self-correcting for nonsingular matrices. It may also
be used to compute the Moore-Penrose pseudoinverse for singular
matrices. However, for such matrices, and for matrices $A(\epsilon)$ where
singular values below $\epsilon$ are ignored, it may be mildly
unstable. The correction step in lines 5-6 of \cref{alg:ns} insures
stability and need only be activated after the unwanted singular values of $XA$ are small enough, a condition that can be inexpensively monitored \cite{pan91}.

\subsection{Newton-Schulz iterations with hierarchical matrices}
To adapt iterative methods to the hierarchical matrix representation context, two modifications are needed as shown in \cref{alg:hns}. The first comes from the fact that the primary fast operation we have, especially on GPUs, is a blocked matrix-vector multiplication; the construction of the Schulz iterates is done via the procedure of \cref{sec:hara} where a sampler, i.e., a matrix-vector expression evaluator, is provided to the construction procedure at every iteration. For example, a simple sampler for NS would follow \cref{alg:hns_sampler} to produce the samples $Y$ from random input vectors.

The second, more consequential modification, is that the iterative
algorithm may be carried out with truncation---an option made possible
by the tunable-accuracy nature of the hierarchical
representation.
As solving exactly the tangent system in the early phases of the Newton method is a waste of computational resources,
it makes little sense to construct the matrices in the first iterations to the ultimately desired accuracy.
In fact, Hackbusch et al.~\cite{hackbusch08} showed that, under
fairly general conditions, the intermediate iterates $X_k$ may be
replaced by approximations without affecting the convergence rate of
the method. Therefore we have the freedom to choose the accuracy
$\epsilon_k$ to which the intermediate iterates are to be
generated. A more effective strategy starts with a low accuracy in the early
iterations and gradually reduces it as convergence is
approached.
As we show below, the use of
such a dynamic threshold reduces substantially the overall
computational cost of inversion by producing intermediate iterates
with smaller footprints.

\begin{algorithm}[!ht]
\caption{Hierarchical matrix iterative inversion of $A$ starting with $X_0$}
\label{alg:hns}
\begin{algorithmic}[1]
\Procedure{hnewton\_schulz}{$A$, $X$, $X_0$, $\epsilon$}
	\State $k = 0,\ X_k = X_0$
	\While{$E = \norm{AX-I}_2 > \epsilon$}
		\State $S = \textsc{ns\_eval}<\! X_k, A \!>$  \Comment{\emph{sampler for (k+1)\textsuperscript{st} iterate}}
		\State $\epsilon_k = \texttt{setContructionThreshold}\left(E, k, \epsilon \right)$ \Comment{\emph{dynamic threshold}}
		\State $X_{k+1} = \texttt{buildH}\left( S, \epsilon_k \right)$  \Comment{\emph{construction of \cref{sec:hara}}}
		\State $k = k + 1$
	\EndWhile
\EndProcedure
\end{algorithmic}
\end{algorithm}

\begin{algorithm}[!ht]
\caption{Templated sampler for Newton-Schulz iterate, computes $Y = X_{k+1} \Omega$}
\label{alg:hns_sampler}
\begin{algorithmic}[1]
\Procedure{ns\_eval<$X_k$,$A$>}{$\Omega$, $Y$}
	\State $Y = (2I - X_k A) X_k \Omega$
\EndProcedure
\end{algorithmic}
\end{algorithm}

We show the effectiveness of \cref{alg:hns} on the following minimal
surface
problem.
Let $m(x)$ be the height of a surface defined in the unit square $\Omega = [0, 1]^2$.
\begin{equation}
\begin{array}{lll}
        \underset{m(x)}{\mbox{minimize }}  &  {\displaystyle J(m) = \int_\Omega \sqrt{ 1 + |\nabla m|^2}} d\Omega, \vspace*{0.0em} \\
        \mbox{subject to} &  m(\partial \Omega) = m_0.
\end{array}
\label{eq:shape_opt_problem}
\end{equation}
We discretize $m(x)$ and the objective using finite differences on an
$n = 128 \times 128$ grid. The Hessian for this problem can be
computed exactly, up to discretization errors, since it corresponds to
the linearization of a nonlinear Poisson equation.
Its inverse, however, is dense, and therefore it is a good test for the effectiveness of the iterative inversion method described above, since all resulting approximation errors are attributable to it.

\begin{figure}[!ht]
	\begin{subfigure}[b]{0.325\textwidth}
		\includegraphics[width=\textwidth]{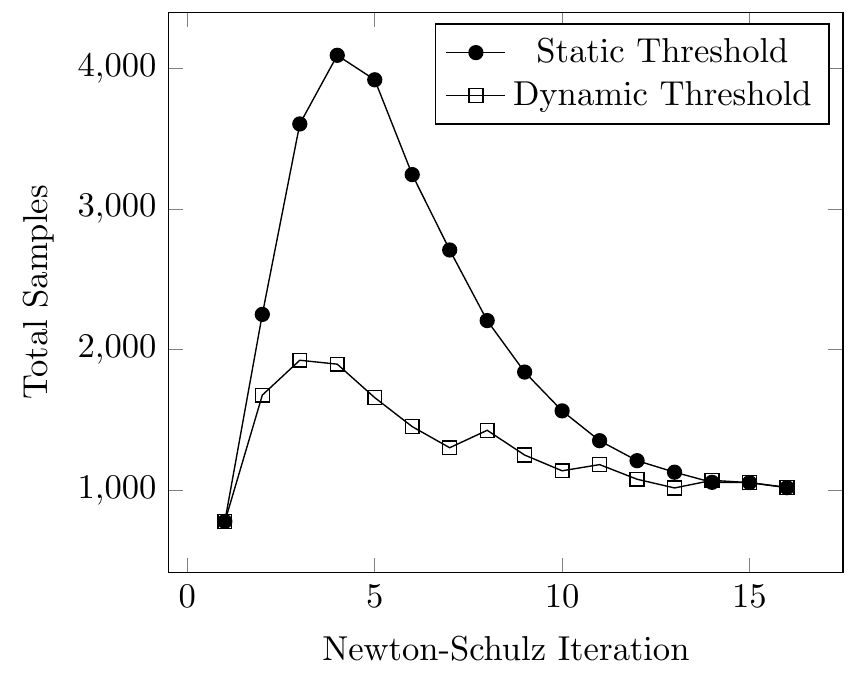}
		\caption{}
		\label{fig:shape_opt_ns_samples}
	\end{subfigure}
	\hfill
	\begin{subfigure}[b]{0.325\textwidth}
		\includegraphics[width=\textwidth]{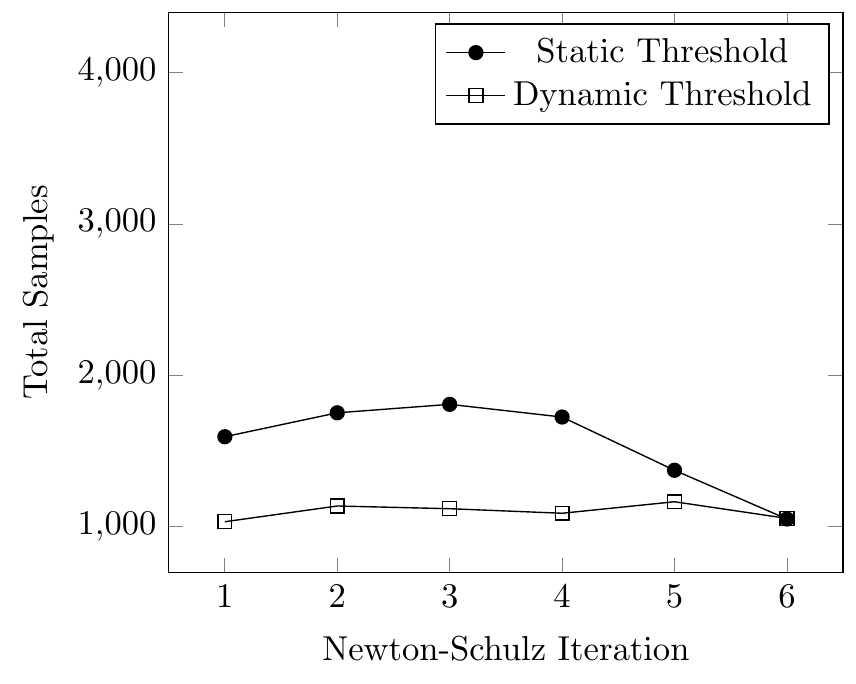}
		\caption{}
		\label{fig:shape_opt_ns_samples_it2}
	\end{subfigure}
	\hfill
	\begin{subfigure}[b]{0.325\textwidth}
		\includegraphics[width=\textwidth]{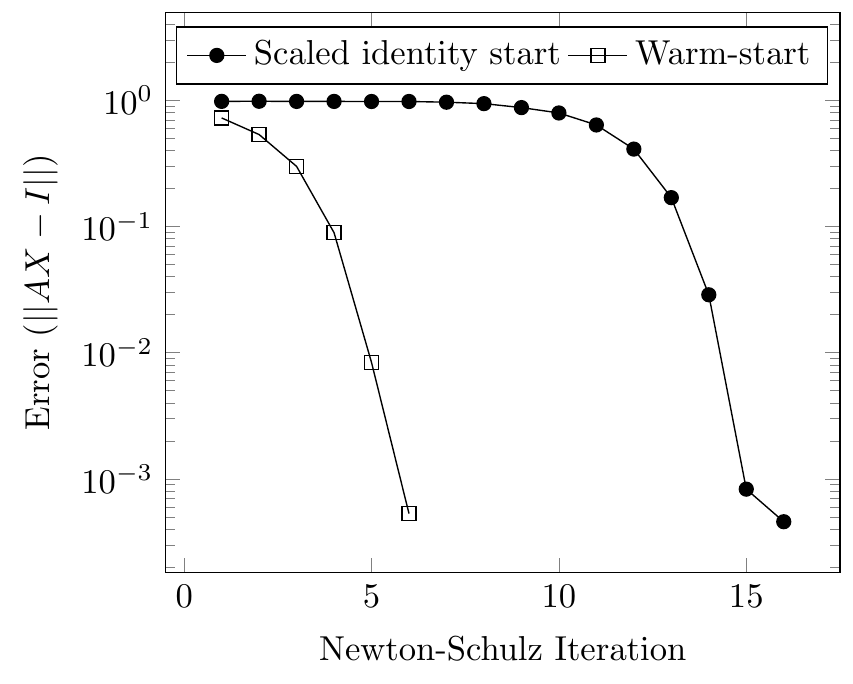}
		\caption{}
		\label{fig:shape_opt_ns_err_it}
	\end{subfigure}
\caption{Performance of Hessian inversion on a problem of size $n = 16{,}384$. Number of Hessian-vector products throughout the NS iterations: (a) starting from a scaled identity and (b) a warm starting from previous iteration.
(c) Convergence history.}
\end{figure}

For the following experiments, we set the compression threshold used during the construction of the hierarchical matrix to $\epsilon_k = 10^{-6}$ with a leaf size of $64$.
Starting from the scaled identity $X_0 = \frac{I}{\norm{A}_\infty}$,
inverting the Hessian using NS for the first optimization step needs
quite a few iterations, with the number of required samples for
intermediate iterates increasing very rapidly before receding as it
converges as shown in \cref{fig:shape_opt_ns_samples}. The total time
needed to invert the Hessian using this method is about 178 seconds on
a P100 GPU, with over $75\%$ of the total runtime spent in
compression. To alleviate
the impact of the intermediate iterates, we can start with a
relatively loose compression error threshold $\epsilon_k$, tightening
the threshold as we converge. This significantly reduces the number of samples needed for the earlier iterations as shown in \cref{fig:shape_opt_ns_samples}, where we start with a much looser threshold of $10^{-2}$. The total runtime is then reduced to 51s, a $3.5\times$ reduction in inversion time, with $69\%$ of the time spent in compression.

From the second optimization iteration onwards, we can use the
approximate Hessian inverse from the previous iteration as an initial
guess, allowing NS to converge in fewer iterations. While the first
optimization iteration converged in 16 NS iterations, the second
iteration needed only 6, inverting the Hessian in 32s if a static
threshold is used and in 14s with a dynamic threshold. In
\cref{fig:shape_opt_ns_samples_it2}, we compare the number of needed samples for static and dynamic thresholds.
Finally, the error in the inverse $\norm{AX-I}_2$ throughout the NS procedure is shown in \cref{fig:shape_opt_ns_err_it} for a given optimization step, exhibiting rapid convergence as we approach the solution. Warm-starting from previous iterations results in the iterations entering the fast convergence regime earlier.

\subsection{Higher order methods for faster convergence}
The results above show that the hierarchical matrix compression
portion of the construction algorithm dominates the total runtime of
the inversion algorithm. This is because we apply the low rank updates generated during the sampling phase in relatively small blocks due to memory constraints of the GPU. The majority of the resulting linear algebra operations during compression, such as batched rank-revealing QR and tall skinny QR decompositions, are not particularly efficient on GPUs and are relatively costly. On the other hand, the sampling operations, which primarily consist of batched matrix-matrix products and blocked sparse matrix vector products, are highly parallel, efficient, and arithmetically intensive. Taking this disparity into consideration, we employ high order hyperpower iterative methods to shift the computational load to the sampling phase. An order $l$ hyperpower iteration is defined as \cite{pan18}
\begin{equation}
X_{k+1} = X_k \left(I + R_k + \dots + R_k^{l-1}\right) = X_k \sum_{i=0}^{l-1} R_k^i,
\label{eq:hyperpower_iteration}
\end{equation}
where $R_k = I - A X_k$, involving $l$ matrix products. Setting $l=2$ gives us the standard NS iteration of the previous section. Most previous work on efficient hyperpower iterations with dense matrices attempt to reduce the number of matrix-matrix products by calculating a few temporary matrices and factoring the summation in \eqref{eq:hyperpower_iteration}. Here, we seek to avoid the costly compression for those temporary matrices, and use the original form of the equation which performs $l$ products. This achieves our goal of concentrating the workload on the far more efficient sampling phase of the computation, and achieves considerable time savings. Sampling $X_{k+1}$ can be done efficiently using a method similar to Horner's method for polynomial evaluation, as shown in \cref{alg:hyperpower_eval} which replaces the NS evaluation of line 4 of \cref{alg:hns}.
\begin{algorithm}[th]
\caption{Sampling hyperpower iterate $X_{k+1}$ of order $l$, computes $Y = X_{k+1} \Omega$.}
\label{alg:hyperpower_eval}
\begin{algorithmic}[1]
\Procedure{hyperpower\_eval<$X_k$,$A$>}{$\Omega$, $Y$, $l$}
	\State $Y = R = \Omega$
	\For{$i = 1 \dots l-1$}
		\State $R = (I - A X_k) R$
		\State $Y = Y + R$
	\EndFor
	\State $Y = X_k Y$
\EndProcedure
\end{algorithmic}
\end{algorithm}

\begin{figure}[!ht]
\begin{subfigure}[b]{0.45\textwidth}
	\includegraphics[width=\textwidth]{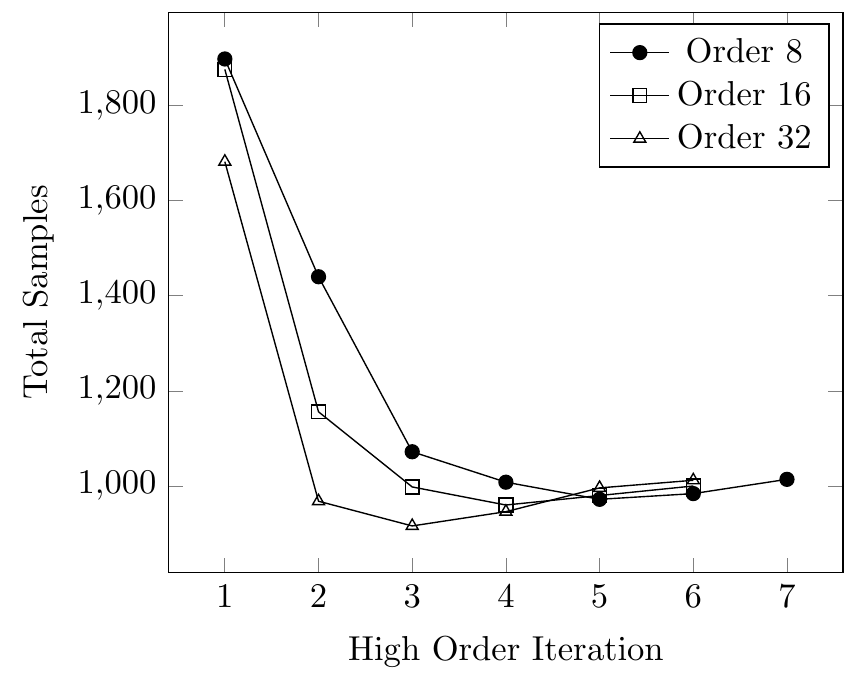}
	\caption{Number of samples vs iterations.}
	\label{fig:high_order_samples}
\end{subfigure}
\hfill
\begin{subfigure}[b]{0.45\textwidth}
	\includegraphics[width=\textwidth]{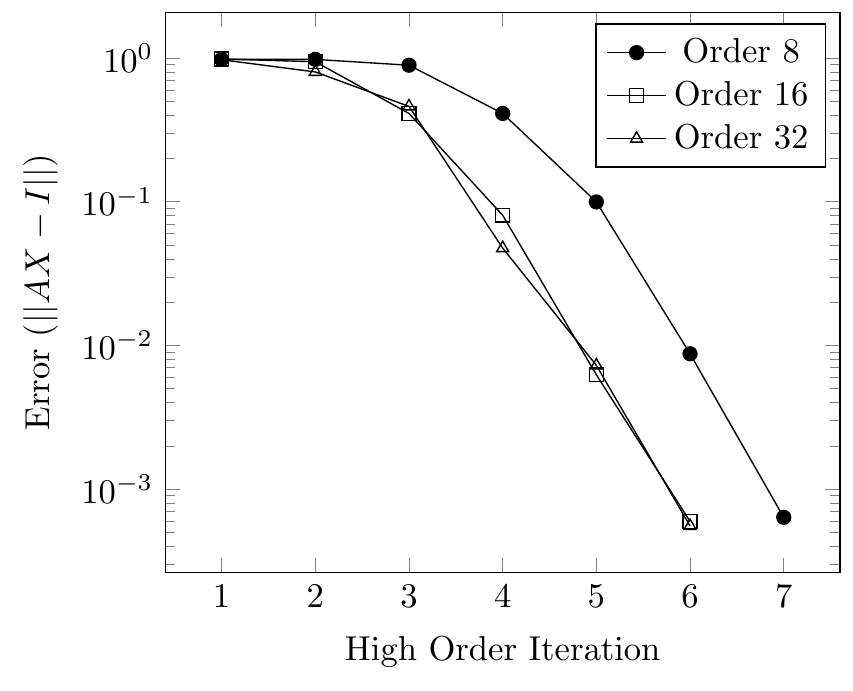}
	\caption{Convergence history.}
	\label{fig:high_order_err}
\end{subfigure}
\caption{Performance of Hessian inversion using high order iterations of order 8, 16, and 32.}
\end{figure}

The performance of higher order methods is illustrated in \cref{fig:high_order_samples} which shows the number of iterations and the samples taken in each iteration for hyperpower iterations of order $l = 8, 16, 32$, which are notably lower than those of \cref{fig:shape_opt_ns_samples}.
High order methods also have the benefit of faster convergence as shown in \cref{fig:high_order_err} where the order 8 method takes 7 iterations to converge and the order 16 and 32 methods take 6 iterations, as opposed to the 16 iterations needed by the second order method.  The overall inversion times are also considerably lower, with order 8 and 16 at 20 seconds and the order 32 at 25s. The order 8 method does one more iteration than the order 16 method, but the lower sampling cost puts it on equal footing, whereas the order 32 performs the same number of iterations as the order 16 method while having higher cost per sample.

\subsection{Unrolling iterations}
As another algorithmic optimization for the example problem \eqref{eq:shape_opt_problem}, we consider the effect of loop unrolling.
The low number of NS iterations needed for the second optimization iteration onwards, due to the use of the previous approximate Hessian inverse as the initial NS iterate, presents another opportunity for greater inversion performance. By unrolling the iterates to express $X_k$ in terms of $X_0$, we can achieve the same goal as the higher order methods, where the workload is shifted to the sampling phase and the increase in samples for the intermediate iterates can be avoided; the $k$-th unrolled iterate can be evaluated as the hyperpower iterates as:
\begin{equation}
\begin{split}
X_k & = X_{k-1} \left(2I + AX_{k-1} \right) \\
	& = X_{k-2} \left(2I + AX_{k-2} \right) \left(2I + AX_{k-2} \left(2I + AX_{k-2} \right) \right)
	 = \dots \\
	& = X_0 \sum_{i=0}^{2^k} {{2^k}\choose{i}} (-A X_0)^i .
\end{split}
\label{eq:unrolled_iteration}
\end{equation}
Since the number of entries in the sum grows exponentially, we can
efficiently unroll only a small number of iterations (say 5 or 6); however, since the number of iterations needed after the first optimization iteration is small, we can effectively unroll all of the required iterations, and obtain the approximate inverse in a single construction. This  further reduces the inversion time from  14 seconds to 4 seconds, giving us another $3.5\times$ improvement on the dynamic error threshold inversion.

We remark that the inversion methods introduced in this section demonstrate the performance benefits of concentrating the computational effort in the arithmetically intensive sampling phase, not only in terms of reducing inversion time as in the unrolled iterations, but also in terms of reducing the total number of iterations required for convergence, as in the higher order methods. Related methods for computing square roots and inverse square roots may be similarly formulated and optimized.


\section*{Acknowledgments}
This work was supported by the King Abdullah University of Science and Technology (KAUST) Office of Sponsored Research (OSR) under Award No: OSR-2018-CARF-3666.

\bibliographystyle{siamplain}
\bibliography{optim,maxwell,ccgo,helmholtz}

\end{document}